\documentclass[11pt]{article}
\usepackage[a4paper, margin=0.8in]{geometry}

\usepackage{authblk}
\usepackage{algorithm}
\usepackage{algorithmic}
\usepackage{amsmath}
\usepackage{amssymb}
\usepackage{hyperref}
\usepackage{rotating}
\usepackage{xcolor}
\usepackage{graphicx}
\usepackage{natbib}
\DeclareMathOperator{\diff}{d}

\usepackage{listings}

\definecolor{DarkBlue}{rgb}{0.00,0.00,0.55}
\definecolor{DarkRed}{rgb}{0.55,0.00,0.00}
\definecolor{DarkGreen}{rgb}{0.00,0.55,0.00}
\definecolor{Purple}{rgb}{0.5, 0.0, 0.5}
\definecolor{BurntOrange}{RGB}{204,85,0}
\definecolor{backcolour}{rgb}{0.95,0.95,0.92}
\definecolor{Plum}{rgb}{0.572, 0.149, 0.561}

\lstdefinelanguage[firedrake]{python}[]{python}{%
   emph=[1]{
      Function,TestFunction,TrialFunction,
      LinearVariationalProblem,LinearVariationalSolver,
      dx,assign,solve,action,derivative,
   },
   emph=[2]{
      subject,time_derivative,transport,diffusion,source_label,linearisation, 
      Term,keep,drop,replace_subject,has_label,label_map, 
   },
   emph=[3]{with,as,None,True},
}

\lstdefinestyle{framed}{frame=tblr}

\usepackage{fancyhdr} 
\fancypagestyle{plain}{
\fancyhf{}
\lfoot{\textcopyright \ Crown Copyright, Met Office}
\cfoot{\thepage}
}
\newcommand\tophline{\hline\noalign{\vspace{1mm}}}
\newcommand\middlehline{\noalign{\vspace{1mm}}\hline\noalign{\vspace{1mm}}}
\newcommand\bottomhline{\noalign{\vspace{1mm}}\hline}

\title{Gusto: a toolkit for compatible finite element dynamical cores}

\author[1,*]{Jemma Shipton}
\author[2]{Thomas M. Bendall}
\author[3]{Colin J. Cotter}
\author[4]{David M. Acreman}
\author[5, 1]{Timothy C. Andrews}
\author[6, 3]{Jack Betteridge}
\author[1, 2]{Alex Brown}
\author[7]{Thomas H. Gibson}
\author[1]{Nell Hartney}
\author[3]{David A. Ham}
\author[3]{Joshua Hope-Collins}
\author[8]{Jeremy A. Pike}
\author[3]{Connor J. Ward}
\author[1]{Dan Witt}

\affil[1]{Department of Mathematics and Statistics, University of Exeter, Exeter, UK}
\affil[2]{Dynamics Research, Met Office, Exeter, UK}
\affil[3]{Department of Mathematics, Imperial College London, London SW7 2AZ, UK}
\affil[4]{Hewlett Packard Enterprise, Lake View, Brierly Furlong, Bristol, BS34 8SW}
\affil[5]{Department of Climate and Space Sciences and Engineering, University of Michigan, Ann Arbor, Michigan}
\affil[6]{Devito Codes, London, UK}
\affil[7]{Advanced Micro Devices, Inc}
\affil[8]{Research Software and Analytics Group, University of Exeter, Exeter, EX4 4RN, UK}

\begin{document}

\maketitle

\begin{abstract}

Gusto is a tool for rapid prototyping of new algorithms based on
compatible finite element methods for discretising the fluid equations
in the dynamical core of weather and climate models. In this context,
the development of new algorithms typically follows a path that
involves testing on a hierarchy of simpler equation sets sharing
important features of the mathematical structure of the equation set
solved by the dynamical core. Gusto provides compatible finite element
spatial discretisations of the hierarchy of equation sets, along with
a range of timestepping schemes. The code is modular so that different
spatial and temporal discretisations can be combined and the same
schemes can be applied to any of the equation sets. This facilitates
fair comparison of the schemes and makes the extension of any
algorithmic developments to the more complex equation sets more
straightforward. In this paper we present Gusto to the community by
describing the algorithms implemented in Gusto and showcasing their
performance through a variety of test problems commonly used in
numerical weather prediction.
\end{abstract}


\section{Motivation}

Compatible finite element methods provide a route to stable,
consistent discretisations of the partial differential equations
governing geophysical fluids \citep{cotter2023compatible}. They
possess many of the desirable properties described by
\citet{staniforth2012horizontal} as essential for the accuracy of the
dynamical core of weather and climate models. This includes conservation of
mass (and potentially other quantities such as energy and enstrophy),
correct wave propagation speeds and avoidance of spurious
computational modes, even on non-orthogonal grids. Obtaining these
properties on a non-orthogonal grid is important because weather and
climate models achieve optimal performance on supercomputers by
performing calculations concurrently, or in parallel, and orthogonal
grids of a spherical manifold will necessarily have singular points
that lead to a degradation in parallel scaling. An example of this is
the commonly-used latitude-longitude grid where the singular points
are located at the poles.


While the required properties of the discretisation guide the choice
of finite element spaces, there are many options for both the spatial
and temporal discretisation and it is not always clear which will be
best for any particular problem. To address this challenge, we present
Gusto: a code library for compatible finite element discretisations
applied to geophysical fluid dynamics equations relevant to numerical
weather and climate prediction. Gusto is built using the Firedrake
finite element software which provides the underlying finite element
machinery. Firedrake itself makes extensive use of the PETSc linear
and nonlinear solver library, which gives Gusto access to
state-of-the-art efficient and scalable solvers. Together this means
that Gusto provides a flexible, extensible, easy-to-use toolkit of
options for rapid prototyping of different spatial and temporal
discretisations. In this paper we introduce Gusto to the community and
summarise its current capabilities.

\subsection{Background}


The finite element method represents fields (e.g. velocity, pressure,
density or temperature) on a mesh made of non-overlapping, conforming
cells, or elements. Within each element, the field is represented as a
sum of coefficients multiplied by polynomial basis functions. The
finite element space is specified by the degree of the polynomial
basis functions and the continuity of these functions between the
elements. The partial differential equations used in weather and
climate modelling have multiple prognostic fields and it is not
necessary, or indeed desirable, to use the same finite element
representation for each field. Using different finite element
representations for each field gives what is known as a mixed finite
element discretisation.

`Compatible' or `mimetic' finite element discretisations are
constructed by choosing the different finite element representations
for each field so that discrete versions of continuous vector calculus
identities hold. These methods are useful in a variety of applications
due to their ability to capture the underlying geometric structure of
the governing equations. In the context of geophysical fluid dynamics,
these methods are useful because they represent geostrophic balance
exactly; this is crucial for reproducing accurate simulations of the
large-scale flow \citep{staniforth2012horizontal}. In addition to
this, there is flexibility within the compatible finite element
framework to choose representations that also give other advantages
such as energy conservation, local mass conservation and discrete
dispersion properties that lead to no spurious pressure modes.

\begin{figure}
  \includegraphics{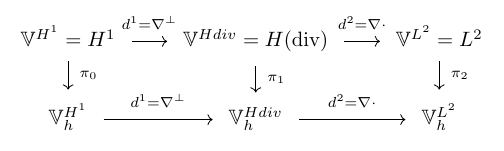}
  \caption{De Rham complex for two dimensional spaces. The
    interelement continuity of the functions decreases from left to
    right, with functions in $\mathbb{V}_h^{H^1}$ being fully
    continuous, those in $\mathbb{V}_h^{Hdiv}$ having continuous normal
    components and those in $\mathbb{V}_h^{L^2}$ being fully
  discontinuous. The polynomial order also decreases from left to
  right.}
  \label{fig:2DdeRham}
\end{figure}

\begin{figure}
  \includegraphics{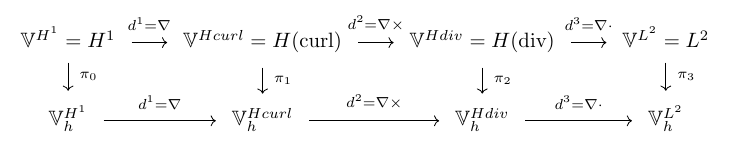}
  \caption{De Rham complex for three dimensional spaces. The
    interelement continuity of the functions decreases from left to
    right, with functions in $\mathbb{V}_h^{H^1}$ being fully
    continuous, those in $\mathbb{V}_h^{Hcurl}$ having continuous
    tangential components, those in $\mathbb{V}_h^{Hdiv}$ having
    continuous normal components and those in $\mathbb{V}_h^{L^2}$ being
    fully discontinuous. The polynomial order also decreases from left
    to right.}
  \label{fig:3DdeRham}
\end{figure}

The compatible finite element spaces relevant for geophysical fluid
dynamics are constructed from a discrete differential complex called
the de Rham complex. This is a sequence of spaces with differential
operators between them that map onto the kernel of the next operator
in the sequence. This is often summarised in a commutation diagram as
in Figures \ref{fig:2DdeRham} and \ref{fig:3DdeRham}. In these
diagrams the infinite dimensional spaces are along the top row and the
$\pi_i$ denote projections into the corresponding finite dimensional
subspace used in the compatible finite element discretisation. The
differential operators are denoted by $\diff^k$ and these commute with
the projections, meaning that $\diff^{k+1}\pi_k(\cdot) =
\pi_{k+1}\diff^{k+1}(\cdot)$. In 3D, this guarantees that $\nabla$
maps onto the kernel of $\nabla\times$, ensuring
$\nabla\times\nabla(\cdot)=0$, and $\nabla\times$ maps onto the kernel
of $\nabla\cdot$, ensuring that $\nabla\cdot\nabla\times(\cdot)=0$;
i.e. the spaces are compatible. In 2D we have that $\nabla^\perp$
(the 2D curl operator) maps onto the kernel of $\nabla\cdot$,
ensuring that $\nabla\cdot\nabla^\perp(\cdot)=0$. It is possible to
write this de Rham complex such that instead $\nabla$ maps onto the
kernel of $\nabla^\perp$, ensuring that
$\nabla^\perp\cdot\nabla(\cdot)=0$, but the former is more useful in
constructing discretisations for geophysical fluid dynamics since it
enables schemes with local mass conservation \citep{cotter2023compatible}.

Due to these promising properties, the next Met Office dynamical core,
Gungho, will use the lowest order compatible finite element
discretisation on the cubed sphere \citep{melvin2024mixed}. However,
there are many important avenues for research and unanswered questions
relating to, for example, the choice of grid and finite element family
for the $H(\textrm{div})$ or $H(\textrm{curl})$ space; the design of
transport schemes; stable and efficient timestepping algorithms;
linear solver techniques and how to couple unresolved physics
processes to the resolved dynamics. This type of research is difficult
in operational models, which are generally designed to execute
specific configurations very efficiently. The purpose of Gusto is to
facilitate this research by providing open source, flexible,
extensible software that provides the key components required when
designing algorithms for a dynamical core in a framework that allows
the same discretisations to be applied to different equation sets from
the hierarchy commonly used to test numerical weather prediction
models. In particular, Gusto will automatically provide the correct
compatible finite element function spaces given the cell shape of the
mesh, the degree of $\mathbb{V}_h^{L^2}$ and the family of
$\mathbb{V}_h^{Hdiv}$. These function spaces can be difficult to
construct on spherical domains where the Piola transform is required
to ensure that vector fields such as the velocity remain tangential to
the surface \citep{rognes2013automating}. An additional complication
for geophysical flows is that the domain is orders of magnitude larger
in the horizontal direction than in the vertical. Combined with the
fact that the atmosphere is close to a state of hydrostatic balance,
this means that meshes must be aligned in vertical columns to ensure
that horizontal errors in the pressure gradient do not contaminate the
representation of hydrostatic balance
\citep{natale2016compatible}. When constructing 3D function spaces in
Gusto (or 2D spaces on a vertical slice domain), we use the vertically
extruded mesh functionality of Firedrake \citep{McRae2016, Bercea2016}
to create the spaces as a tensor product of the horizontal and
vertical spaces. This means that we can easily construct a function
space for the thermodynamic variable that has degrees of freedom
co-located with the vertical velocity; the finite element analogue of
the Charney-Phillips staggering known to work well in finite
difference numerical weather prediction models
\citep{guerra2016high,melvin2018choice}. In addition to setting up the
compatible finite element function spaces, Gusto provides stable and
accurate transport and diffusion schemes for variables with different
interelement continuity (see Section \ref{sec:transport}), along with
a selection of physics parametrisation schemes and a range of time
discretisations (see Section \ref{sec:results} for a range of
examples). Much of the functionality described above is enabled by the
Firedrake finite element library (see Section \ref{sec:design} for
specific details). The close integration of Firedrake, and therefore
Gusto, with the PETSc library of linear and nonlinear solvers
facilitates easy use of these efficient and scalable solvers for any
of the systems of equations arising from the discretisation.

Gusto has grown out of code developed in proof-of-concept studies that
investigated the properties of the compatible finite element framework
for different equation sets. Originally, \citet{cotter2012mixed}
showed that compatible finite element methods preserve the advantages
of the C-grid finite difference discretisation of the linear shallow
water equations on non-orthogonal grids; namely that energy and mass
are conserved, there are no spurious pressure modes and geostrophic
modes are steady. Various formulations of the nonlinear shallow water
equations have been explored, with \citet{mcrae2014energy,
  bauer2018energy, wimmer2020energy} presenting energy and enstrophy
conserving schemes and \citet{shipton2018higher} demonstrating a
potential vorticity conserving scheme. Energy-vorticity conserving
schemes were presented in \citet{wimmer2021energy}, whilst an
energy-entropy conserving scheme was presented in
\citet{ricardo2024entropy}. The default shallow water discretisation
used in Gusto is closer to the discretisation of the 3D equations used
by Gungho that does not require computation of a vorticity
variable. This shallow water discretisation is comprehensively
documented in \citet{hartney2025exploring}. Successfully extending a
discretisation from the 2D shallow water equations to equations
representing a vertically stratified fluid requires accurate
representation of hydrostatic balance. \citet{natale2016compatible}
proved that compatible finite element discretisations on vertically
aligned meshes result in a well defined hydrostatic pressure, meaning
that the discretisation will preserve hydrostatic balance and not be
polluted by spurious pressure modes. The properties of the 3D
discretisation are often first demonstrated using a 2D vertical slice
domain. Some results in this domain using the transport schemes
described in Section \ref{sec:transport} are presented in
\citet{natale2016compatible}; \citet{yamazaki2017vertical} applied the
method to the Eady model, demonstrating its ability to capture
developing fronts; and \citet{cotter2023compatible2} presented a
comprehensive range of vertical slice simulations using a fully
implicit timestepping scheme. \citet{bendall2019recovered} designed a
second order advection scheme for use with lowest order elements,
where standard finite element schemes would only be first
order. \citet{bendall2020compatible} extended the 3D discretisation to
include moist processes.

Bringing together these methods in a well-tested and documented code
base that is modular and can be easily extended provides a robust and
useful tool for investigating the many open research questions in the
design of accurate and efficient methods within the compatible finite
element framework. A number of recent works have used Gusto for
developing and implementing new algorithms and demonstrate the utility
of this approach. \citet{hartney2026exploringa} extended the compatible
finite element shallow water model to include both thermal and moist
processes in order to explore different formulations of the equations
as well as options for coupling these new physics processes to the
dynamical equations
\citep{hartney2026exploringb}. \citet{witt2026splitting} exploited the
tensor product structure of the function spaces to design methods that
have a different finite element order in the horizontal and
vertical. \citet{andrews2026conservative} developed a conservative
transport scheme for scalar tracers that ensures consistency and
\citet{brown2025fast} demonstrated the accuracy of spectral deferred
correction timestepping schemes for fast-wave slow-wave splittings of
the compressible Euler equations.

In addition to the numerical discretisation research questions
described above, there is a rapidly developing field of research in
combining numerical models with data driven algorithms. Recent
developments show that it may be advantageous to couple the numerical
model to the machine learning algorithm during the training phase,
which requires the numerical model to be differentiable
\citep{kochkov2024neural}. Although operational weather and climate models
often have an adjoint or tangent linear model implementation, the
implementation is specific to the configuration of the forward model,
making it difficult to couple to machine learning algorithms. However,
as Gusto is built on top of Firedrake, it inherits the
\textit{automatic} differentiability of this library without any
changes to the code, making it an ideal software framework in which to
pursue research in this area. Firedrake provides capabilities to
couple to external machine learning libraries
\citep{bouziani2023physics, bouziani2024differentiable,
  farsi2026missing}. Automatic differentiability also facilitates
adjoint sensitivity analysis and variational data assimilation. All of
these topics are active areas of research and will be reported
elsewhere.

In this article we will describe the current capabilities of Gusto and
the design features that enable our flexible implementation whereby
different spatial and temporal discretisations can be applied to a
range of different equation sets relevant for weather and climate
modelling. In Section \ref{sec:equations} we will outline three
different equation sets commonly used for modelling geophysical fluids
in order to demonstrate the similarities between the equation sets and
justify the use of a unified codebase. The advantages of Gusto's
unified approach becomes clear in Sections \ref{sec:FEM} and
\ref{sec:timestepping} where we present the compatible finite element
discretisation of the different types of terms (e.g. pressure
gradient, transport, diffusion) seen in Section \ref{sec:equations}
and the generic timestepping implementation. In Section \ref{sec:design} we describe the code structure of Gusto followed by examples
of the Form Manipulation Language (FML) that facilitates symbolic
manipulation of the finite element forms within our timestepping
code. Section \ref{sec:results} contains a selection of numerical
results that showcase the range of discretisations available within
Gusto. We summarise our current status and future directions in
Section \ref{sec:summary}.

\section{Governing Equations}
\label{sec:equations}
The aim of this section is to introduce the three governing equation sets
commonly used in the development of discretisation methods for weather
and climate models, followed in Section \ref{sec:FEM} by enough
detail about the finite element discretisations provided by Gusto to
demonstrate the similarities between them. These similarities guide
the design of the code structure described in the next section and
motivate use of the Form Manipulation Language (FML) that underlies
the flexibility of Gusto, examples of which will be given in Section
\ref{sec:design}. The convention in Gusto is to write the equations
in residual form, meaning that all the terms appear on the left
side. We will follow this convention below in contrast to many papers
which present the equations with the pressure gradient term and any
terms relating to parametrised physics processes on the right side.

\subsection{Shallow Water Equations}
The shallow water equations are commonly used for testing numerical
algorithms and for exploring geophysical fluid dynamics concepts under
simplified conditions. They are useful in both contexts since they are
the simplest set of equations that model motion on both the slow,
geostrophic timescale and the much faster timescale of the
inertia-gravity waves. This timescale separation poses challenges to
the numerical discretisation and also provides a rich enough range of
dynamics to explain many of the phenomena observed in the atmosphere
and oceans \citep{zeitlin2018geophysical}.

The shallow water equations describe the flow of a shallow layer of
fluid subject to gravitational and, optionally, rotational forces. The
prognostic variables are the velocity $\boldsymbol{u}$ 
and the fluid depth $D = H + h - B$ where $H$ is the
undisturbed depth of the fluid layer, $h$ is the free surface
deviation from $H$, and $B$ is the bottom topography, if present. The
equations are:

\begin{align}
  \frac{\partial \boldsymbol{u}}{\partial t} + (\boldsymbol{u}\cdot\nabla)\boldsymbol{u} + f\hat{\boldsymbol{k}}\times\boldsymbol{u} + g\nabla (D+B) &= 0, \\
  \frac{\partial D}{\partial t} + \nabla\cdot(D\boldsymbol{u}) &= 0, \label{eqn:sw_D}
\end{align}
where $f$ is the Coriolis parameter, $\hat{\boldsymbol{k}}$ is the unit vector
pointing upwards and $g$ is the acceleration due to gravity.

\subsection{Compressible Boussinesq Equations}
The compressible Boussinesq equations are commonly used in ocean
modelling. They are useful from a numerical discretisation perspective
because they incorporate compressibility effects without the
complexity of discretising the nonlinear pressure gradient terms that
we shall see in the compressible Euler equations. We write these
equations in terms of the prognostic variables $\boldsymbol{u}$, $p$ and $b$
where $\boldsymbol{u}$ is the velocity field, $p$ is pressure and $b$ is the
fluid buoyancy. In terms of these variables, the equations are:

\begin{align}
  \frac{\partial \boldsymbol{u}}{\partial t} +
  (\boldsymbol{u}\cdot\nabla)\boldsymbol{u} +
  2\boldsymbol{\Omega}\times \boldsymbol{u} + \nabla p - b\hat{\boldsymbol{k}} & = 0, \\
  \frac{\partial p}{\partial t} + c^2\nabla\cdot\boldsymbol{u} & = 0, \\
  \frac{\partial b}{\partial t} + (\boldsymbol{u}\cdot\nabla) b & = 0,  \label{eqn:bouss_b}
\end{align}
where $\boldsymbol{\Omega}$ is the rotation vector of the domain,
$\hat{\boldsymbol{k}}$ is the unit vector pointing upwards and $c$ is the
acoustic wave speed. The equations describe the motion of a
three dimensional stratified fluid but can also be solved in a
two dimensional vertical slice domain described by $(x, z)$
coordinates. In this case, the velocity $\boldsymbol{u}$ can either be
two dimensional as well (if there is no Coriolis force), or fully
three dimensional but with each Cartesian component only varying in
the $x-z$ plane. In the latter case, this is implemented using a
three dimensional computational domain that is only one element wide
in the $y$-direction and with periodic boundary conditions also in
that direction.

\subsection{Compressible Euler Equations}
The dynamical core at the heart of any weather or climate model solves
the rotating compressible Euler equations that model atmospheric
flow. We write these equations in terms of the prognostic variables
$\boldsymbol{u}$, $\rho$ and $\theta$ where $\boldsymbol{u}$ is the velocity field,
$\rho$ is the fluid density and $\theta$ is the potential
temperature. The pressure gradient term is written in terms of the
Exner pressure $\Pi$ which is calculated through an equation of
state. In terms of these variables, the equations are:

\begin{align}
  \frac{\partial \boldsymbol{u}}{\partial t} + 
  (\boldsymbol{u}\cdot\nabla)\boldsymbol{u} +
  2\boldsymbol{\Omega}\times \boldsymbol{u} + c_p\theta\nabla \Pi + g\hat{\boldsymbol{k}} & = 0, \\
  \frac{\partial \rho}{\partial t} + \nabla\cdot(\boldsymbol{u}\rho) & = 0,  \label{eqn:Euler_rho}\\
  \frac{\partial \theta}{\partial t} + (\boldsymbol{u}\cdot\nabla)\theta & = 0,  \label{eqn:Euler_theta}\\
  \Pi^{(1-\kappa)/\kappa} & = \frac{R}{p_0}\rho\theta, & 
\end{align}
where $\Pi$ is the Exner pressure, $p_0$ is a reference pressure, $R$
is the gas constant, $c_p$ is the specific heat at constant pressure
and $\kappa=R/c_p$. Again, these equations describe the motion of a
three dimensional stratified fluid but can also be solved in a
two dimensional vertical slice domain described by $(x, z)$
coordinates as described above.

In all of the equation sets described here, the nonlinear velocity
advection term $(\boldsymbol{u}\cdot\nabla)\boldsymbol{u}$ can be replaced by the sum
of a vorticity-based term and the gradient of the kinetic energy. This
is known as the vector invariant form and it enables the construction
of schemes that conserve energy and/or enstrophy
\citep{mcrae2014energy, bauer2018energy, wimmer2020energy,
  wimmer2021energy}. However, this form can exhibit numerical
instabilities \citep{bell2017numerical} but as there are persuasive
arguments for each form, in Gusto we provide both forms of the equations.

\subsection{Physics parametrisations}
\label{sec:physics}
The equations described so far in this section govern the adiabatic
geophysical fluid motion on scales that are resolved by the spatial
and temporal discretisation. Weather and climate models also include
physics parametrisations representing unresolved, diabatic
processes. For example, \citet{bendall2020compatible} describes the
extension of the compressible Euler equations to include moist
thermodynamics and \citet{hartney2026exploringa} explores a variety of
ways to represent moist thermodynamics in the shallow water
equations. The simplest example of a physics parametrisation is the
addition of a diffusive or viscous term of the form $\nu\nabla^2q$ for
any prognostic variable $q$ to any of the corresponding equations
described above, where $\nu$ is the diffusion or viscosity
coefficient. We will describe how to discretise terms of this form in
the following section.

\section{Finite Element Discretisation}
\label{sec:FEM}

The discretisations in Gusto use the method of lines, meaning that the
spatial and temporal dimensions are decoupled and can be discretised
separately. In this section we consider the spatial
discretisation. Comparing the equation sets given in the previous
section, we see several similarities. For example, the prognostic
equation for the velocity always contains a nonlinear transport term,
a term due to rotation and a `pressure' gradient term (in shallow
water this is the gradient of the depth but the term plays the same
role in the dynamics and is derived from the pressure gradient term in
the full equations \citep{zeitlin2018geophysical}). The other
prognostic fields satisfy either an advective or continuity form of
the transport equation. These similarities remain when we formulate
the compatible finite element discretisations, as we shall see below,
and importantly this means that the discretisation of these terms can
be written and implemented in software so that it can be applied to
any of the equation sets. In the following, we refer to terms arising
from the transport, or advection, operator
$\boldsymbol{u}\cdot\nabla q$ or $\nabla\cdot(\boldsymbol{u}q)$ for any field $q$ as `transport' terms; terms
arising from the pressure gradient force or Coriolis force as
`forcing' terms; and any terms arising from parametrisations of
unresolved physics processes as `physics' terms.

The finite element discretisation is formed by choosing appropriate
representations, or function spaces, for the prognostic variables,
multiplying by test functions from these function spaces and
integrating over the domain. In addition to this, we integrate by
parts when the interelement continuity of the basis functions is not
sufficient to define the result of applying the differential operator
directly. For the domains considered in this paper, the only boundary
conditions required are those of no normal flow, $\boldsymbol{u}\cdot\boldsymbol{n}=0$
where $\boldsymbol{n}$ is the outward pointing normal to the boundary. This
means that the test function for the velocity must also be chosen from
a function space that has this boundary condition. We always represent the domain using Cartesian coordinates in physical space so that the vertical coordinate is the physical height, not a pressure-based coordinate. The formulation of the equations is separate from the mesh geometry but our standard choice for constructing meshes that have a vertical coordinate is to extrude the horizontal mesh. This results in a mesh that has vertically-aligned columns, which improves representation of hydrostatic balance, and is structured for optimal implementation of the line smoothers required as part of the multigrid component of the iterative solvers. When the domain is spherical, or an extruded spherical shell, Firedrake produces a piecewise polynomial approximation of the sphere. In this paper, we use a quadratic representation of the coordinate field in order to achieve the correct convergence rates. Using higher order finite element spaces for the prognostic variables would require a higher order representation of the coordinates.

In line with the theory outlined in the introduction, we always choose
the velocity $\boldsymbol{u} \in \mathring{\mathbb{V}}_h^{Hdiv}$ with
$\mathring{\mathbb{V}}_h^{Hdiv} \subset \mathbb{V}_h^{Hdiv} \subset
H(\text{div})$ and, as mentioned above, the $\mathring{\phantom{\mathbb{V}}}$ indicates that
the functions in $\mathring{\mathbb{V}}_h^{Hdiv}$ satisfy the no normal flow
boundary conditions (if required). The field whose gradient appears in the
`pressure' gradient term is chosen to be in the space
$\mathbb{V}_h^{L^2}$ with $\mathbb{V}_h^{L^2} \subset L_2$, or, in the
case of the more complex nonlinear pressure gradient term for the
compressible Euler equations, it is a function of a variable in this
space. The compressible Boussinesq and compressible Euler equations
have a third prognostic variable, either buoyancy or potential
temperature, which is chosen to be in the space $\mathbb{V}_h^\theta$
which has degrees of freedom that are co-located with those of the
vertical velocity, mimicking the Charney-Phillips staggering as
discussed in \citet{melvin2018choice}. These choices are governed by
the dispersion properties of the linearised equation sets. To solve the
nonlinear equation sets we need to introduce stable and accurate
discretisations of the transport terms and the nonlinear pressure
gradient term in the compressible Euler equation set. In addition,
some test cases require viscous and diffusive terms to generate
converged solutions (for example \citet{straka1993numerical}). In the
following subsections we illustrate a non-exhaustive selection of the
options available in Gusto.

\subsection{Pressure gradient term}
\label{sec:pressuregradient}
The pressure gradient term is discretised in a similar way for each
equation set, with some additional considerations required for the
nonlinear form of this term in the compressible Euler equations. Since
the pressure (or depth, for shallow water) field is discretised in
$\mathbb{V}_h^{L^2} \subset L_2$, meaning that it is discontinuous
between cells, the gradient is not defined globally and the term must
be integrated by parts. In the case that this term is linear, using
the notation for the compressible Boussinesq equations, the weak form
is

\begin{equation}
\int_\Omega\boldsymbol{w}\cdot\nabla p \diff x = -\int_\Omega\nabla\cdot\boldsymbol{w}p\diff x, \quad \forall \boldsymbol{w} \in \mathring{\mathbb{V}}^{Hdiv},
\end{equation}
where $\Omega$ denotes the domain and there are no surface terms from
the integration by parts because these vanish on interior element
boundaries (since the normal component of $\boldsymbol{w}$ is continuous
across element facets) and on exterior boundaries we set no normal
flow boundary conditions. However, the pressure gradient term in the
compressible Euler equations, $c_p\theta\nabla\Pi$, is nonlinear and
this introduces an additional complication. Integration by parts is
still required since $\Pi$ is evaluated as a function of the density
$\rho \in \mathbb{V}_h^{L^2}$ and $\theta \in \mathbb{V}_h^\theta$ so is
discontinuous between elements, but then integrating by parts
transfers the differential operator to $\boldsymbol{w}\theta$ which is
discontinuous in the horizontal direction. Following
\citet{natale2016compatible} we integrate by parts in each element
$e$, giving

\begin{equation}
  c_p\int_e\boldsymbol{w}\cdot\theta\nabla\Pi\diff x = -c_p\int_e\nabla\cdot(\boldsymbol{w}\theta)\Pi\diff x + c_p\int_{\partial e}\hat{\theta}\boldsymbol{w}\cdot\boldsymbol{n}\hat{\Pi}\diff S, \quad \forall \boldsymbol{w} \in \mathring{\mathbb{V}}^{Hdiv},
\end{equation}
where $\partial e$ denotes the element facets, $\boldsymbol{n}$ is the outward
pointing normal to the facet and $\hat{q}$ denotes the value of a
variable $q$ on the facet, which for discontinuous variables does not
have a unique definition. Summing over the domain gives

\begin{equation}
  c_p\int_\Omega\boldsymbol{w}\cdot\theta\nabla\Pi\diff x = -c_p\int_\Omega\nabla\cdot(\boldsymbol{w}\theta)\Pi\diff x + c_p\int_{\Gamma_v}[\![\theta\boldsymbol{w}\cdot\boldsymbol{n}]\!]\{\Pi\} \diff S, \quad \forall \boldsymbol{w} \in \mathring{\mathbb{V}}^{Hdiv},
\end{equation}
where $\Gamma_v$ denotes the vertical element facets, across which
$\theta$, and therefore $\Pi$, are discontinuous. Here we have chosen,
as in \citet{natale2016compatible}, to take the average value of $\Pi$
across the facets, $\{\Pi\} = (\Pi^++\Pi^-)/2$, where the $+$ and
$-$ indicate each side of the facet. The brackets $[\![\boldsymbol{v}]\!]=\boldsymbol{v}^+\cdot\boldsymbol{n}^++\boldsymbol{v}^-\cdot\boldsymbol{n}^-$
denote the jump in a vector quantity $\boldsymbol{v}$ across a facet, where $+$ and $-$ label each side of the facet, with $\boldsymbol{n}^\pm$ being the unit normal pointing out of the corresponding element. The integral over the
horizontal facets vanishes since both $\theta$ and $\boldsymbol{w}\cdot\boldsymbol{n}$
are continuous across these facets so the jump is zero.

\subsection{Transport terms}
\label{sec:transport}
The transport terms are those of the form $(\boldsymbol{u}\cdot\nabla)q$ or
$\nabla\cdot(\boldsymbol{u}q)$ for some prognostic variable $q$. Since the
prognostic variables are discretised using different finite element
function spaces $\mathbb{V}_h^{Hdiv}$, $\mathbb{V}_h^{L^2}$ and
$\mathbb{V}_h^\theta$, each with different continuity properties
between the elements, the choice of stable transport discretisation
can be different for each prognostic variable. The spatial order of
the finite element spaces also imposes requirements on the choice of
transport scheme. In particular, when the lowest order space
$\mathbb{V}_h^{L^2}$ is chosen to be piecewise constant, we need to
use the recovered transport scheme described in
\citet{bendall2019recovered} to achieve second order convergence
overall. Here we describe some transport schemes used in Gusto for the
next-to-lowest order sequence of function spaces (i.e. where the
lowest order space $\mathbb{V}_h^{L^2}$ is chosen to be piecewise
linear and discontinuous between elements) and refer the reader to
\citet{bendall2019recovered, bendall2023improving} for details of the
recovered schemes, which are also implemented in Gusto.

\subsubsection{Transport of $\mathbb{V}_h^{L^2}$ fields}
\label{subsec:DGtransport}
The shallow water and compressible Euler equation sets both contain a
prognostic equation that is the continuity form of the transport
equation for a variable represented by a function that is fully
discontinuous between elements (see equation~(\ref{eqn:sw_D}) and equation~(\ref{eqn:Euler_rho}) respectively). This means that the divergence term is
not globally defined. A stable discretisation of this equation can be
constructed by integrating that term by parts and taking the upwind
values on the element facets. Using the compressible Euler equation
notation, this gives, in each element $e$ with boundary $\partial e$:

\begin{equation}
  \int_e \phi\nabla\cdot(\rho\boldsymbol{u}) \diff x = -\int_e \nabla\phi \cdot (\rho\boldsymbol{u}) \diff x + \int_{\partial e} \phi\hat{\rho}\boldsymbol{u}\cdot\boldsymbol{n} \diff S, \quad \forall \phi \in \mathbb{V}_h^{L^2},
\end{equation}
where $\hat{\rho}$ denotes the value of $\rho$ on the element
boundary. Summing over the domain and chosing
$\hat{\rho}$ to be the upwind value of $\rho$ for stability gives

\begin{equation} \label{eqn:conservative_DG_upwind}
  \int_\Omega \phi\nabla\cdot(\rho\boldsymbol{u}) \diff x = -\int_\Omega \nabla\phi \cdot (\rho\boldsymbol{u}) \diff x + \int_\Gamma [\![\phi\boldsymbol{u}]\!]\tilde{\rho} \diff S, \quad \forall \phi \in \mathbb{V}_h^{L^2},
\end{equation}
where $\Omega$ is the domain, $\Gamma$ denotes the element facets,
$\tilde{\rho}$ denotes the value of $\rho$ on the upwind side of the
facet and the brackets $[\![\boldsymbol{v}]\!]$
denote the jump in a vector quantity $\boldsymbol{v}$ across a facet as defined in section \ref{sec:pressuregradient}. \\
\\
However, when solving equation~(\ref{eqn:conservative_DG_upwind}) using a multi-stage
timestepping scheme to evolve a density field in time, the total increment
over a timestep will include the nested evaluations of differential operators
upon the velocity.
As a result, a density field that is initially constant in space will not
receive an increment equal to $-\rho\nabla\cdot\boldsymbol{u}$.
The importance of this property is discussed in \citet{bendall2025swift} as
`divergence linearity', and is related to the desirable `divergence-consistency'
property of conservative Semi-Lagrangian schemes discussed by \citet{wood2026replicating}.
This problem can be addressed by using an `advective-then-flux' formulation,
in which the advective form of the transport equation is used for all but the final
stages of a multi-stage scheme to avoid taking the divergence of a density which
already includes an update from the velocity divergence.
The advective form is discretised as:
\begin{equation} \label{eqn:advective_DG_upwind}
  \int_\Omega \phi\boldsymbol{u}\cdot\nabla\rho \diff x = -\int_\Omega \rho\nabla\cdot(\phi\boldsymbol{u}) \diff x + \int_\Gamma \phi \boldsymbol{u}\cdot\boldsymbol{n}\tilde{\rho} \diff S, \quad \forall \phi \in \mathbb{V}_h^{L^2}.
\end{equation}
The final stage of the timestepping scheme then uses the conservative (also known as `flux') form of the
equation~(\ref{eqn:conservative_DG_upwind}), using the density $\rho$ which is a linear sum of
terms from the previous advective stages.
This ensures that whole scheme is conservative, but also obeys the `divergence linearity'
property.
This approach has been implemented in Gusto, and for more details see the descriptions
in \citet{bendall2023solution} and \citet{melvin2024mixed} in which this `advective-then-flux'
form was used with finite volume transport schemes.

\subsubsection{Transport of $\mathbb{V}_h^\theta$ fields}
\label{subsec:thetatransport}
The compressible Boussinesq and compressible Euler equation sets both
contain a prognostic equation that is the advective form of the
transport equation for a variable that is represented by a function
that is continuous in the vertical direction but discontinuous in the
horizontal direction (see equation~(\ref{eqn:bouss_b}) and equation~(\ref{eqn:Euler_theta}) respectively). This means that the gradient of the field
required in the transport term is not globally defined. There are two
common options implemented in Gusto for discretising this term.
The first is the `embedded DG' scheme presented by \citet{cotter2016embedded}, which
constructs the transported $\mathbb{V}_h^\theta$ field in a discontinuous (or broken)
form of the $\mathbb{V}_h^\theta$ space, before using the upwind DG transport
scheme as in equation~(\ref{eqn:advective_DG_upwind}).
\\
\\
As described in \citet{yamazaki2017vertical}, a second stable discretisation of
this equation can be constructed by using a combination of the above
upwinding strategy in the discontinuous horizontal direction and a
streamline-upwind Petrov-Galerkin (SUPG) method
\citep{brooks1982streamline} in the continuous vertical direction. The
SUPG method replaces the finite element test function $\gamma$ with
$\gamma + \tau \boldsymbol{u}\cdot\nabla\gamma$ where $\tau$ is a
stabilisation parameter, although here we just require the vertical
part of the gradient term. This has the effect of suppressing spurious
oscillations associated with the continuous finite element
approximation of the transport term. The stabilisation parameter
$\tau$ is set to be $\Delta t/\sqrt{15}$ by default in Gusto,
following \citet{raymond1976selective}.

We first integrate by parts in each column $C$ of the mesh, taking the
upwind value of the transported field $\tilde{\theta}$ on the vertical
facets of the mesh, denoted by $\Gamma_v$, giving

\begin{align}
  \int_C \gamma\boldsymbol{u}\cdot\nabla\theta\diff x &= -\int_C \nabla\cdot(\gamma\boldsymbol{u})\theta \diff x + \int_{\Gamma_v}\gamma \boldsymbol{u}\cdot\boldsymbol{n}\tilde{\theta}\diff S, \quad \forall \gamma \in \mathbb{V}_h^\theta, \\
  &= \int_C \gamma\boldsymbol{u}\cdot\nabla\theta \diff x + \int_{\Gamma_v}\gamma \boldsymbol{u}\cdot\boldsymbol{n}(\tilde{\theta}-\theta_{int})\diff S, \quad \forall \gamma \in \mathbb{V}_h^\theta,
\end{align}
where we have integrated by parts a second time to move the
differential operator from $\gamma$ back to $\theta$ since we need to
replace the test function $\gamma$ with $\gamma+\tau w\gamma_z$, where
$w=\boldsymbol{u}\cdot\hat{\boldsymbol{z}}$ is the vertical component of the
velocity. This gives a second boundary term that now contains the
value of $\theta$ in the interior of the column, denoted by
$\theta_{int}$. Summing over the domain and applying the SUPG method
in the vertical gives
\begin{equation}
  \int_\Omega \gamma\boldsymbol{u}\cdot\nabla\theta\diff x =\int_\Omega(\gamma+\tau w\gamma_z)\boldsymbol{u}\cdot\nabla\theta \diff x + \int_{\Gamma_v}[\![(\gamma+\tau w\gamma_z) \boldsymbol{u}]\!]\tilde{\theta}-[\![(\gamma+\tau w\gamma_z) \boldsymbol{u}\theta]\!]\diff S, \quad \forall \gamma \in \mathbb{V}_h^\theta,
\end{equation}
where the brackets $[\![\boldsymbol{v}]\!]$ denote the jump in a vector
quantity $\boldsymbol{v}$ across a facet as defined in section \ref{sec:pressuregradient}. Note that, as explained in
\citet{natale2017scale} and \citet{cotter2023compatible}, the test
function should also be replaced in the time derivative term for
consistency.

\subsubsection{Transport of $\mathbb{V}_h^{Hdiv}$ fields}
All three of the equation sets described in Section \ref{sec:equations} include a prognostic equation for the velocity field
$\boldsymbol{u}\in\mathbb{V}_h^{Hdiv}$ that contains a nonlinear transport term
$\boldsymbol{u}\cdot\nabla\boldsymbol{u}$. This term can be discretised directly using
a vector version of the scheme in Section \ref{subsec:thetatransport} or, as in \citet{natale2018variational} and
\citet{gibson2019compatible}, it can be replaced by the vector
invariant form and then discretised. In Gusto we offer both options
and we present the vector invariant scheme below.

The starting point for this discretisation is to replace the
$\boldsymbol{u}\cdot\nabla\boldsymbol{u}$ term with $(\nabla\times\boldsymbol{u})\times\boldsymbol{u}
+ \frac{1}{2}\nabla(\boldsymbol{u}\cdot\boldsymbol{u})$. Now, multiplying by a test
function and integrating by parts over an element $e$ with boundary
$\partial e$ gives
\begin{multline}
  \int_e \boldsymbol{w}\cdot(\nabla\times\boldsymbol{u})\times\boldsymbol{u}
  + \frac{1}{2}\boldsymbol{w}\cdot\nabla(\boldsymbol{u}\cdot\boldsymbol{u}) \diff x \\= \int_e\nabla\times(\boldsymbol{u}\times\boldsymbol{w})\cdot\boldsymbol{u}\diff x - \int_{\partial e}\boldsymbol{n}\times(\boldsymbol{u}\times\boldsymbol{w})\cdot\boldsymbol{u} \diff S - \frac{1}{2}\int_e(\nabla\cdot\boldsymbol{w})(\boldsymbol{u}\cdot\boldsymbol{u})\diff x, \quad \forall \boldsymbol{w} \in \mathbb{V}_h^{Hdiv}.
\end{multline}
Summing over the domain and taking the upwind value of $\boldsymbol{u}$,
denoted by $\tilde{\boldsymbol{u}}$, on the element boundaries gives
\begin{multline}
  \int_\Omega \boldsymbol{w}\cdot(\nabla\times\boldsymbol{u})\times\boldsymbol{u}
  + \frac{1}{2}\boldsymbol{w}\cdot\nabla(\boldsymbol{u}\cdot\boldsymbol{u}) \diff x \\= \int_\Omega\nabla\times(\boldsymbol{u}\times\boldsymbol{w})\cdot\boldsymbol{u}\diff x - \int_{\Gamma}[\![\boldsymbol{n}\times(\boldsymbol{u}\times\boldsymbol{w})]\!]\cdot\tilde{\boldsymbol{u}} \diff S - \frac{1}{2}\int_\Omega(\nabla\cdot\boldsymbol{w})(\boldsymbol{u}\cdot\boldsymbol{u})\diff x, \quad \forall \boldsymbol{w} \in \mathbb{V}_h^{Hdiv},
\end{multline}
where the brackets $[\![\boldsymbol{v}]\!]$ denote the jump in a vector
quantity $\boldsymbol{v}$ across a facet as defined in section \ref{sec:pressuregradient}.

\subsubsection{Summary of transport discretisation}
The transport schemes described above do not preserve
monotonicity unlike the continuous scalar transport equations~(\ref{eqn:sw_D}), (\ref{eqn:bouss_b}), and
(\ref{eqn:Euler_rho})-(\ref{eqn:Euler_theta}), so the discretisation can
generate spurious overshoots and undershoots. When monotonicity is
important, for example, to avoid generating spurious and unphysical
negative values of moisture variables, we can apply a limiter as part
of the transport scheme. These limiters are described in
\citet{bendall2020compatible}.

It is common in weather and climate models for the horizontal and
vertical directions to be discretised differently due to the different
scales of motion in each direction. Typically the vertical part of the
transport term is treated implicitly in the time stepping scheme so
that the timestep is not restricted by the vertical Courant number
(e.g. as in \citet{weller2013runge}). We provide the option to split
the transport terms into horizontal and vertical components, as
described in \citet{brown2025fast}.

\subsection{Viscous and diffusive terms}
\label{sec:diffusion}
For the simulations presented here we do not require viscous or
diffusive terms for stability, but sometimes they are required to
produce converged solutions so that refining the grid does not lead to
smaller and smaller scale motions. These terms also provide a simple example of a physics scheme as described in section \ref{sec:physics}. Viscous or diffusive terms involve
the Laplacian of the prognostic field, which is not
globally defined for all the finite element function spaces we use in this
paper. In Gusto we use the symmetric interior penalty discretisation
method \citep{cockburn2007note} which we describe below applied to the
viscous term added to the prognostic velocity equation when
required. This method can be adapted to apply to scalar fields.

First we multiply the term $\nu\nabla^2\boldsymbol{u}$ by a test function
$\boldsymbol{w}\in\mathbb{V}_h^{Hdiv}$, where $\nu$ is the viscosity
coefficient, assumed to be constant. Integrating by parts introduces
surface terms which, after applying the symmetric interior penalty
method \citep{cockburn2007note} become:

\begin{equation}
 \nu\int_\Omega\nabla\boldsymbol{w}:\nabla\boldsymbol{u}\diff x - \int_\Gamma[\![\boldsymbol{v}\otimes\boldsymbol{n}]\!]:\{\nabla\boldsymbol{u}\}\diff S - \int_\Gamma[\![\boldsymbol{u}\otimes\boldsymbol{n}]\!]:\{\nabla\boldsymbol{v}\}\diff S + \mu\int_\Gamma [\![\boldsymbol{u}\otimes\boldsymbol{n}]\!]:[\![\boldsymbol{v}\otimes\boldsymbol{n}]\!]\diff S,
\end{equation}
where $:$ denotes double contraction, $\{q\}$ and $[\![\boldsymbol{v}]\!]$
denote the average and jump across a facet (both defined in section \ref{sec:pressuregradient}), and $\mu$ is the penalty parameter,
typically chosen to be $10/h$ where $h$ is a measure of the gridsize.

\section{Timestepping}
\label{sec:timestepping}

We provide a range of standard timestepping schemes in Gusto,
including explicit, implicit and implicit-explicit Runge-Kutta schemes
and multilevel schemes such as Adams-Moulton and Adams-Bashforth
schemes. The modular structure of the Gusto software facilitates the design and implementation of other timestepping schemes such as the 
spectral deferred correction methods described in
\citet{brown2025fast}. However, in general, we are particularly
motivated by the use of a Semi-Implicit Quasi-Newton (SIQN) timestepping
method that closely follows the scheme used at the Met Office in both the
ENDGame and Gungho dynamical cores \citep{wood2014inherently, melvin2019mixed,
  melvin2024mixed, johnson2026regional} as this is stable over longer
timesteps. This algorithm is described in detail in those references
and we provide an overview in the next section.

\subsection{Semi-Implicit Quasi-Newton timestepping scheme} \label{sec:siqn}

This scheme has evolved from a tradition of Semi-Implicit
Semi-Lagrangian (SISL) schemes, which have been used successfully for operational
numerical weather forecasting for several decades.
SISL schemes combine a Crank-Nicolson time
discretisation with Semi-Lagrangian transport to facilitate a large
stable timestep with minimal computational effort. The motivation behind this approach is to enable use of accurate explicit timestepping schemes for the nonlinear transport terms while treating the stiff terms arising from the forcing implicitly. Consider the general
equation, i.e. the weak form of any of the equation sets introduced in
Section \ref{sec:equations}, in the form

\begin{equation}
  \frac{\partial \chi}{\partial t} +  \mathcal{T}(\chi) + \mathcal{F}(\chi) + \mathcal{P}(\chi) = 0,
\end{equation}
where $\chi$ is the vector of prognostic variables, while $\mathcal{T}$, $\mathcal{F}$ and $\mathcal{P}$ denote the transport, forcing and physics terms respectively. A number of options for coupling physics parametrisations $\mathcal{P}$  to the SIQN timestepping
have been implemented in Gusto. While a full description can be found in
\citet{hartney2026exploringb}, here and in Algorithm \ref{alg:SIQN} we focus on the SIQN timestepping for $\mathcal{T}$ and $\mathcal{F}$ and include only
physics parametrisations $\mathcal{P}$ that are evaluated at the end of the timestep.
Following the SISL approach, the value of $\chi$ at the $(n+1)$-th time level
is found from the value at the $n$-th time level according to:
\begin{equation} \label{eqn:sisl}
\chi^{n+1}+\alpha\Delta t\mathcal{F}(\chi^{n+1})
= \mathcal{T}_{\overline{\boldsymbol{u}}}\left[\chi^n - (1-\alpha)\Delta t\mathcal{F}(\chi^n)\right],
\end{equation}
where $\mathcal{T}_{\overline{\boldsymbol{u}}}$ evaluates the transport of prognostic variables
by the velocity $\overline{\boldsymbol{u}}=\tfrac{1}{2}\left(\boldsymbol{u}^n+\boldsymbol{u}^{n+1}\right)$,
and $\alpha$ is the semi-implicit off-centering parameter. A SISL scheme would use a Semi-Lagrangian transport scheme to evaluate $\mathcal{T}_{\overline{\boldsymbol{u}}}$. Instead, Gusto uses Eulerian transport schemes which
can be subcycled to extend the overall timestep.
The Met Office's ENDGame and Gungho dynamical cores solve equation~(\ref{eqn:sisl}) using a
`Quasi-Newton' approach, which iterates towards the solution $\chi^{n+1}$ via
a nested loop structure.
Stages in the timestep are outlined in Algorithm \ref{alg:SIQN}, using
notation defined in this section.
Once any explicit contributions arising from the forcing terms have been evaluated, a timestep consists of an `outer'
loop, indexed here by $i$, in which the transporting velocity is updated with the latest value of
$\mathbf{u}^{n+1}$, and prognostic variables are transported.
Once variables have been transported, the algorithm enters the `inner' loop, indexed here by $k$,
which involves the evaluation of implicit forcing terms to obtain a residual
$\chi_r$, defined as
\begin{equation}
  \chi_r = -\big(\chi^{n+1}_{i,k} + \alpha\Delta t \mathcal{F}(\chi^{n+1}_{i,k}) - \chi_T\big).
\end{equation}
As described by \citet{hartney2026exploringa}, the full Newton method would
iterate towards the solution by solving
\begin{equation}
\mathcal{J}[\chi_{i,k}^{n+1}](\chi_{i,k+1}^{n+1} - \chi_{i,k}^{n+1}) = \chi_r,
\end{equation}
where $\mathcal{J}$ is the Jacobian of $\chi_r$ with respect to $\chi^{n+1}$,
and $k$ denotes the iteration index ($i$ being the `outer' iteration index that is unchanged here).
In the `Quasi-Newton' method, the Jacobian is approximated using a linearisation
$\mathcal{S}$ of the equation set, which gives a linear system to solve in each
inner iteration:
\begin{equation} \label{eqn:linear_system}
\mathcal{S}(\chi_{i,k+1}^{n+1} - \chi_{i,k}^{n+1}) = \mathcal{S}\chi'_{i,k} = \chi_r.
\end{equation}
This `incremental' form follows that used in Gungho, and is
described in more detail by \citet{melvin2019mixed,melvin2024mixed}.
In Gusto, the strategy for solving the linear system
in equation~(\ref{eqn:linear_system}) is to approximately eliminate any variables
in the $\mathbb{V}_h^\theta$ space, and then introduce a hybridisable
formulation for the velocity and $\mathbb{V}_h^{L^2}$ variable by
breaking the velocity continuity constraints and introducing trace
Lagrange multipliers to enforce them in the solution again. We then eliminate down to the trace variables using the SLATE system in
Firedrake \citep{gibson2020slate}. The linear systems used in Gusto are described in full by
\citet{bendall2020compatible,hartney2025exploring}.
The use of Firedrake and PETSc allows a wide range of possible solver techniques
that could be used for solving this linear problem; by default Gusto uses the
generalized minimal residual method (GMRES), preconditioned with a horizontal
multigrid with column line smoothers, as presented in \citet{betteridge2023hybridised}. Appendix \ref{sec:siqn_details} describes some remaining details about the
formulation of the Semi-Implicit Quasi-Newton scheme used in Gusto.

\begin{algorithm}
  \caption{Semi-Implicit Quasi-Newton timestepping scheme.}
  \label{alg:SIQN}
  \begin{algorithmic}[1]
    \STATE Set: $\chi^{n+1}_0 = \chi^{n}$
    \STATE Explicit forcing: $\chi_{f_e} = \chi^n - (1-\alpha)\Delta t \mathcal{F}(\chi^n)$
    \FOR {Outer loop: $i < i_{max}$}
      \STATE Update transporting velocity: $\overline{\boldsymbol{u}}_i = \frac{1}{2} (\boldsymbol{u}_i^{n+1} + \boldsymbol{u}^n)$
      \STATE Explicit transport: $\chi_T = \mathcal{T}_{\overline{\boldsymbol{u}}_i}(\chi_{f_e})$
      \STATE Set: $\chi^{n+1}_{i,0} = \chi^{n+1}_{i}$
      \FOR {Inner loop: $k < k_{max}$}
        \STATE Solve: $\mathcal{S}\chi'_{i,k} = -\big(\chi^{n+1}_{i,k} + \alpha\Delta t \mathcal{F}(\chi^{n+1}_{i,k}) - \chi_T\big)$ for $\chi'_{i,k}$
        \STATE Increment: $\chi^{n+1}_{i,k+1} = \chi^{n+1}_{i,k} + \chi'_{i,k}$
      \ENDFOR
      \STATE Set: $\chi^{n+1}_{i+1} = \chi^{n+1}_{i, k_{max}}$  
    \ENDFOR
    \STATE Explicit physics: $\chi^{n+1} = \mathcal{P}(\chi^{n+1}_{i_{max}})$
    \STATE Advance timestep: $\chi^n = \chi^{n+1}$
  \end{algorithmic}
\end{algorithm}

\section{Software Design}
\label{sec:design}
Gusto is built on top of the Firedrake package, an automated system
for solving partial differential equations using the finite element
method \citep{FiredrakeUserManual}. This gives many advantages, not
least that it provides the finite element machinery required to
translate the discretisations written in section \ref{sec:FEM}
into the linear and nonlinear equations that we need to solve. The
weak forms are specified using the Unified Form Language (UFL), a
domain specific language that enables the definition of the finite
element forms within Python \citep{alnaes2014unified}. Firedrake takes
these forms and generates fast, efficient C code
\citep{Luporini2016} using TSFC \citep{Homolya2017} ensuring
computational efficiency by reordering meshes
\citep{LangeMitchellKnepleyGorman2015}. The other packages that
Firedrake depends on are summarised in \citet{davies2022towards} but
we note that we make particular use of higher order mesh
representation on curved manifolds \citep{Homolya2017a}, quadrilateral
meshes \citep{Homolya2016, McRae2016}, extruded meshes
\citep{McRae2016, Bercea2016} and the PETSc library of scalable linear
and nonlinear solvers \citep{petsc-web-page, petsc-user-ref,
  petsc-efficient} via a Python interface \citep{Dalcin2011},
including use of custom block preconditioning \citep{kirby2018solver}.

\begin{figure}
  \centering{\includegraphics[width=0.8\textwidth]{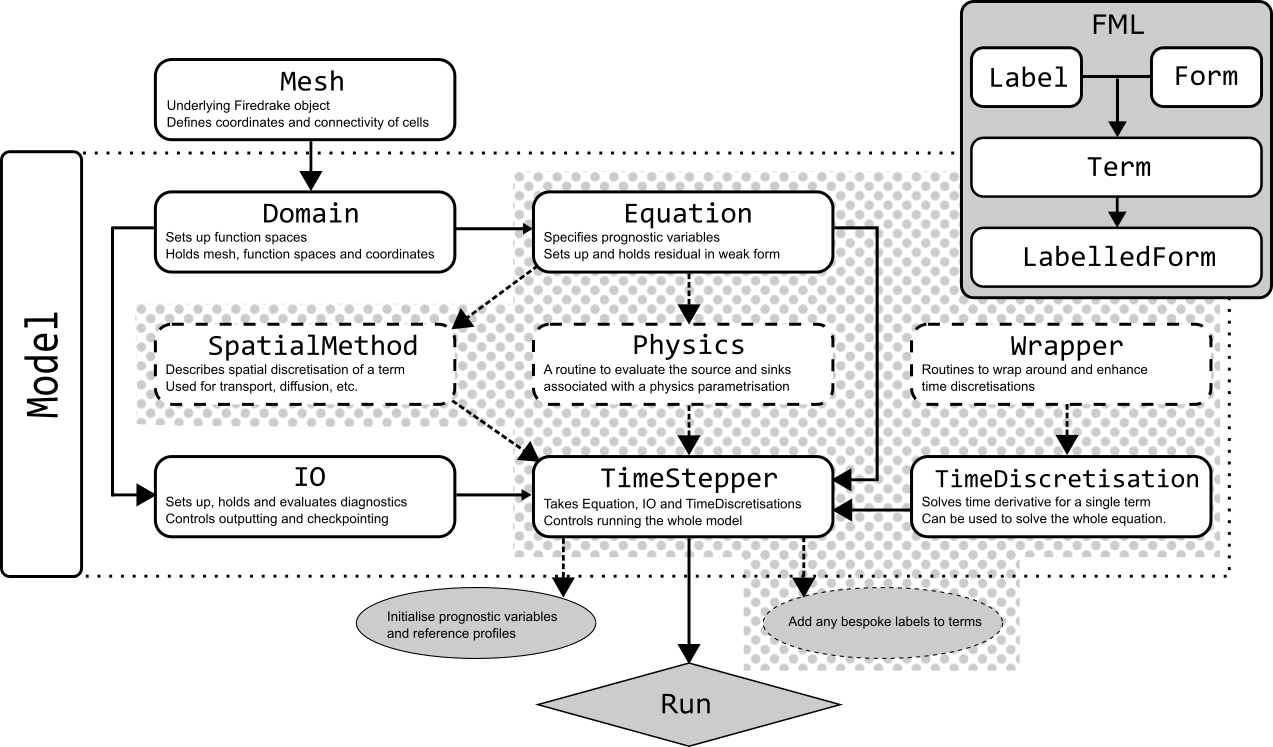}}
  \caption{A schematic summarising the code design of Gusto. Boxes describe
  different types of object defined within the code, while shaded ellipses
  represent steps that the user may define, and the `run' method is highlighted
  in a diamond.
  Boxes with dashed outlines are optional and may not be used in all configurations.
  Arrows indicate the how the objects interact, with the direction indicating that
  an object is an argument passed to another object.
  The `Model' object can encapsulate most of the other objects, which is indicated
  via a dotted box.
  The FML infrastructure and the objects defined within it are shown in the top
  right corner; the dotted background highlights the parts of Gusto that are
  underpinned by FML.}
  \label{fig:code_schematic}
\end{figure}

Here we outline the design of the Gusto code, and the software
building blocks that play different roles within a generic model,
illustrated in Figure \ref{fig:code_schematic}. These building blocks
can be encapsulated within a \texttt{Model} object which allows users
to easily run a range of test problems with the same numerical
discretisation options. The Gusto model comprises a domain and a
prognostic equation set. The \texttt{Domain} class holds the mesh and the finite
element function spaces along with any other information about the
space-time domain such as the timestep. The prognostic equation sets
are each defined in an \texttt{Equation} class, making use of the common finite
element forms for the transport and diffusion terms to set up the
residual (i.e. the left hand side of the equations defined in Section
\ref{sec:equations}) in weak form. At this point, no particular
discretisation of the transport or diffusion forms is specified. This
happens in the \texttt{SpatialMethod} class which can be used to apply
the methods described in sections \ref{sec:transport} and \ref{sec:diffusion} to the relevant terms. These spatial methods are required
to specify a stable discretisation for any nonlinear equation set, but
do not need to be specified for some linearisations of the
equations. Any physics source terms, such as those required for moist
processes or turbulence parametrisation, are added by calling the
\texttt{Physics} class that defines them. \texttt{TimeDiscretisation} classes define schemes for advancing fields in time, such as the Runge-Kutta, Adams-Moulton and Adams-Bashforth schemes. If required, the \texttt{TimeDiscretisation} can be enhanced by a wrapper that implements additional steps, such as mapping to higher order function spaces, as required by the embedded DG scheme \citep{cotter2016embedded}. The \texttt{Timestepper}
class contains the definition of the timestepping procedure and the
\texttt{run} method that loops over the timesteps. It also sets up
diagnostic fields and calls the IO as required. The simplest
timestepping class will just apply one time discretisation method to
the equation set but it is common in weather and climate prediction
models to split the different types of terms and apply a time
discretisation separately to each type of term. For example, in the
Semi-Implicit Quasi-Newton scheme the transport and forcing terms are
solved separately. \\ \\ In Algorithm \ref{alg:basic} we show the
basic structure of the timestepping loop in Gusto. The first step is
to modify the equation residual to apply the spatial discretisation
(for example, the DG upwind or discretisation described in Section
\ref{sec:transport}). Before entering the timestepping loop, we set
up the IO, including any diagnostics, and write out the initial
conditions. In this basic timestepper, the timestepping loop just
applies the time discretisation to the equation set. Any output is
then generated at frequencies specified by the user. We output
diagnostic fields and a time series of global diagnostics, such as
minimum and maximum values of the fields. We also write out
checkpointing files that store fields in a format that can be read back
into Firedrake functions to compute further diagnostics, or to compare
with other simulations (as in a convergence test) or to initialise a
new simulation.

Algorithms \ref{alg:SIQN} and \ref{alg:splitphysics} illustrate two
more complex timestepping algorithms. These replace lines 4-5 in
Algorithm \ref{alg:basic}. In Algorithm \ref{alg:splitphysics} we
separate the dynamics terms from the physics terms, solving first for
the dynamics and then sequentially applying each physics term in the
order in which they were specified by the user. As described in section \ref{sec:siqn}, in
Algorithm \ref{alg:SIQN} the splitting goes further, separating the
dynamics into transport terms and forcing terms.

\begin{algorithm}
  \caption{Basic timestepping scheme}
  \label{alg:basic}
  \begin{algorithmic}[1]
    \STATE Replace transport and diffusion terms with specified spatial discretisation scheme.
    \STATE Setup IO and diagnostics and write out initial conditions.
    \WHILE{$t < t_{max}$}
      \STATE Apply timestepping scheme to equation set: $\chi^{n+1} = \mathsf{T}(\chi^n, \chi^{n+1})$
      \STATE Advance timestep: $\chi^n = \chi^{n+1}$
      \IF{output}
        \STATE Update timeseries diagnostics file.
        \STATE Compute and write out diagnostic fields.
        \STATE Write out checkpoint file.
      \ENDIF
      \STATE Increment the current time: $t \gets t + \Delta t$
    \ENDWHILE
  \end{algorithmic} 
\end{algorithm}

\begin{algorithm} 
  \caption{Split physics timestepping scheme}
  \label{alg:splitphysics}
  \begin{algorithmic}[1]
    \STATE Dynamics step: $\hat{\chi}_0 = \mathcal{D}(\chi^n)$
    \FOR {$\mathcal{P}$ in physics schemes}
      \STATE Apply physics scheme: $\hat{\chi}_{i+1} = \hat{\chi}_{i+1} + \mathcal{P}(\hat{\chi}_{i+1})$
    \ENDFOR
    \STATE End of timestep: $\chi^{n+1} = \hat{\chi}_{np}$ where $np$ is the number of physics schemes
    \STATE Advance timestep: $\chi^n = \chi^{n+1}$
  \end{algorithmic} 
\end{algorithm}

\subsection{Setting up a Gusto model}
In this section we describe how to set up the Gusto model and briefly
demonstrate the flexibility of the software. We begin with code to
solve the rotating shallow water equations on the sphere for the
geostrophic balance test case (Case 2 from
\citet{williamson1992standard}) - see Listing \ref{lst:sw}. To begin
with, we import from both the Firedrake and Gusto modules. Then we set
up the required parameters such as the size of the Earth, mean depth
of the fluid and maximum wind velocity. This example uses an
icosahedral mesh of the sphere with 16 triangular cells along each of
the panel edges and, as required for accuracy, a quadratic
approximation to the coordinate field. The \texttt{Domain} class holds
the mesh and timestep and sets up the compatible finite element
function spaces, given the family (in this case Brezzi-Douglas-Marini)
of the $H(\textrm{div})$ space and the degree (in this case 1) of the
$L^2$ space. The \texttt{OutputParameters} class allows the user to
specify the name of the output directory and the frequency of output
(in timesteps). The \texttt{ShallowWaterParameters} holds the
parameters for this specific simulation, including commonly-used
Earth-like default values for the strength of gravity and the
planetary rotation rate. We now instantiate the
\texttt{ShallowWaterEquations} class which will set up the residual
for this equation set. We specify the trapezium rule as the
timestepper by passing an instance of this class to the
\texttt{Timestepper} class, along with a list of spatial methods for
each prognostic variable. The rest of the code involves the
specification of the initial conditions, in this case a zonal flow
velocity and corresponding geostrophically balanced depth field, as
specified in \citet{williamson1992standard}. The final line runs the
timestepper for 5 days.

\begin{minipage}{0.94\linewidth}
\vspace{0.5em}
\begin{lstlisting}[
   language={[firedrake]python},
   basicstyle=\linespread{0.8}\ttfamily\footnotesize,
   style=framed, frame=single,
   firstline=10,
   caption={A complete script setting up the Gusto model and solving
  the rotating shallow water equations on an Earth-sized sphere.  (See \texttt{sw\_icos\_geostrophic\_balance.py} in \citet{thomas_bendall_2026_24173215})}
]
from firedrake import *
from gusto import *

# Parameters for this case
radius = 6371220                   # planetary radius (m)
mean_depth = 5960                  # reference depth (m)
day = 24*60*60                     # length of day (s)
u_max = 2*pi*radius/(12*day)       # Max amplitude of the zonal wind (m/s)

# Set up model domain
# Use an icosahedral mesh with 16 cells along each panel edge and quadratic (degree 2) coordinates
mesh = GeneralIcosahedralSphereMesh(radius, num_cells_per_edge_of_panel=16, degree=2)
# Specify 900s timestep and BDM family for HDiv space, degree 1 elements for L2 space
domain = Domain(mesh, dt=900, family='BDM', degree=1)

# Set up I/O to write out every 96 timesteps (1 day) to results/sw_geostrophic_balance directory
output = OutputParameters(dirname='sw_geostrophic_balance', dumpfreq=96)
io = IO(domain, output)

# Set up parameters and equation
parameters = ShallowWaterParameters(mesh, H=mean_depth)
eqns = ShallowWaterEquations(domain, parameters)

# Set up timestepper
spatial_methods = [DGUpwind(eqns, 'u'), DGUpwind(eqns, 'D')]
stepper = Timestepper(eqns, TrapeziumRule(domain), io, spatial_methods)

# Velocity and depth fields for initialisation
u0 = stepper.fields('u')
D0 = stepper.fields('D')

# Set parameters and coordinates
g = parameters.g                  # use default value of gravity for Earth
Omega = parameters.Omega          # use default value of planetary rotation rate for Earth
xyz = SpatialCoordinate(mesh)     # Cartesian coordinates of mesh
_, lat, _ = lonlatr_from_xyz(xyz[0], xyz[1], xyz[2])   # spherical coordinates, only need latitude
e_lon = xyz_vector_from_lonlatr(1.0, 0.0, 0.0, xyz)    # unit longitudinal basis vector

# Initialise zonal flow velocity (equations 90-91 of Williamson et al, (1992))
u0.project(u_max*cos(lat)*e_lon)
# Initialise balanced fluid depth (equation 95 of Williamson et al, (1992))
D0.interpolate(mean_depth - (radius * Omega * u_max + 0.5*u_max**2)*(sin(lat))**2/g)

# Run for 5 days
stepper.run(t=0, tmax=5*day)
\end{lstlisting}
\label{lst:sw}
\end{minipage}

The flexibility of Gusto has two key advantages. Firstly, the ease
with which different methods can be applied to the solve the same
equations, meaning that new methods can be directly compared to
existing methods with minimal changes to user code. Secondly, the ease
with which the same method can be applied to solve different
equations, meaning that a newly developed method can quickly progress
through the hierarchy of numerical weather prediction test cases.

To illustrate the first advantage we describe the changes required to
use different spatial and temporal discretisations to solve the
shallow water equations illustrated in Listing \ref{lst:sw}. To
change the spatial discretisation to use Raviart-Thomas elements on
quadrilateral cells requires changing the mesh and the family of
elements passed into the \texttt{Domain} class (see Listing \ref{lst:sw_cs}). To change the temporal discretisation to use the classical
fourth order Runge Kutta method requires changing the scheme passed to
the \texttt{Timestepper} class (see Listing \ref{lst:sw_rk4}). In
each case, all other code from Listing \ref{lst:sw} remains
unchanged.

\begin{minipage}{0.94\linewidth}
\vspace{0.5em}
\begin{lstlisting}[
   language={[firedrake]python},
   basicstyle=\linespread{0.8}\ttfamily\footnotesize,
   style=framed, frame=single,
   linerange={20-23}, firstnumber=11,
   caption={An example showing how to adapt Listing \ref{lst:sw} to use
  Raviart-Thomas elements on quadrilateral cells. (See \texttt{sw\_cs\_geostrophic\_balance.py} in \citet{thomas_bendall_2026_24173215})}
]
# Use a cubed sphere mesh with 16 cells along each panel edge and quadratic (degree 2) coordinates
mesh = GeneralCubedSphereMesh(radius, num_cells_per_edge_of_panel=8, degree=2)
# Specify 900s timestep and RTCF family for HDiv space, degree 1 elements for L2 space
domain = Domain(mesh, dt=900, family='RTCF', degree=1)
\end{lstlisting}
\label{lst:sw_cs}
\end{minipage}

\begin{minipage}{0.94\linewidth}
\vspace{0.5em}
\begin{lstlisting}[
   language={[firedrake]python},
   basicstyle=\linespread{0.8}\ttfamily\footnotesize,
   style=framed, frame=single,
   linerange={33-35}, firstnumber=24,
   caption={An example showing how to adapt Listing \ref{lst:sw} to use the classical fourth order Runge Kutta timestepping method. (See \texttt{sw\_cs\_geostrophic\_balance.py} in \citet{thomas_bendall_2026_24173215})}
]
# Set up timestepper
spatial_methods = [DGUpwind(eqns, 'u'), DGUpwind(eqns, 'D')]
stepper = Timestepper(eqns, RK4(domain), io, spatial_methods)
\end{lstlisting}
\label{lst:sw_rk4}
\end{minipage}

For the SIQN algorithm described in section \ref{sec:siqn}, we provide
a \texttt{Model} class that encapsulates our recommended options (see
Appendix \ref{sec:siqn_details}). Listing \ref{lst:sw_model} shows how to
set up this class for the next-to-lowest order spaces (i.e. the same
spaces as in the original Listing \ref{lst:sw}). Note that we no longer
need to set up the \texttt{Domain} class because \texttt{SIQNModel}
takes care of this for us, and the \texttt{ShallowWaterEquations}
class is not instantiated as this is also taken care of inside the
\texttt{SIQNModel} class.

\begin{minipage}{0.94\linewidth}
\vspace{0.5em}
\begin{lstlisting}[
   language={[firedrake]python},
   basicstyle=\linespread{0.8}\ttfamily\footnotesize,
   style=framed, frame=single,
   linerange={21-33}, firstnumber=11,
   caption={An example showing how to adapt Listing \ref{lst:sw} to use the SIQN algorithm. (See \texttt{sw\_model.py} in \citet{thomas_bendall_2026_24173215})}
]
# Use an icosahedral mesh with 16 cells along each panel edge and quadratic (degree 2) coordinates
mesh = GeneralIcosahedralSphereMesh(radius, num_cells_per_edge_of_panel=16, degree=2)

# Set up I/O to write out every 96 timesteps (1 day) to results/sw_geostrophic_balance directory
output = OutputParameters(dirname='sw_geostrophic_balance', dumpfreq=96)

# Set up parameters and equation
parameters = ShallowWaterParameters(mesh, H=mean_depth)
eqns = ShallowWaterEquations

# Set up model
model = SIQNModel(mesh, dt, parameters, eqns, family='BDM')
model.setup(output)
\end{lstlisting}
\label{lst:sw_model}
\end{minipage}

Now we illustrate the second advantage by adapting Listing \ref{lst:sw_model} to solve the compressible Euler equations on the sphere
(see Listing \ref{lst:compressible_model}). We can use the same base
finite element discretisation but now we need to extrude the mesh in
the vertical direction to generate a 3D domain. The only other change
required is to use the classes defining the parameters and equations.

\begin{minipage}{0.94\linewidth}
\vspace{0.5em}
\begin{lstlisting}[
   language={[firedrake]python},
   basicstyle=\linespread{0.8}\ttfamily\footnotesize,
   style=framed, frame=single,
   linerange={20-33}, firstnumber=20,
   caption={An example showing how to adapt Listing \ref{lst:sw_model} to solve the compressible Euler equations using the SIQN algorithm. (See, for example, \texttt{dry\_baroclinic\_wave.py} in \url{https://github.com/firedrakeproject/gusto_case_studies}) }
]
# Use an icosahedral mesh with 16 cells along each panel edge and quadratic (degree 2) coordinates
base_mesh = GeneralIcosahedralSphereMesh(radius, num_cells_per_edge_of_panel=16, degree=2)
mesh = ExtrudedMesh(base_mesh, layers=10, layer_height=100, extrusion_type='radial')

# Set up I/O to write out every 96 timesteps (1 day) to results/sw_geostrophic_balance directory
output = OutputParameters(dirname='compressible_euler', dumpfreq=96)

# Set up parameters and equation
parameters = CompressibleParameters(mesh)
eqns = CompressibleEulerEquations

# Set up model
model = SIQNModel(mesh, dt, parameters, eqns, family='BDM')
model.setup(output)
\end{lstlisting}
\label{lst:compressible_model}
\end{minipage}

\subsection{Form Manipulation Language}
\label{sec:FML}
In the previous section we demonstrated the ease with which users can
set up simulations in Gusto using different discretisation methods and
equation sets. This flexibility is enabled by a modular code structure
where the equations are defined separately from the timestepping
scheme. This powerful functionality is facilitated through the form
manipulation language (FML) which is used to associate a label to each
term in the residual form defined in the equation class. The
timestepping methods are then able to reason about which manipulations
to apply to which forms. Common manipulations include: 1) filtering
based on the presence of a label or its value, 2) detecting components
of the weak form, such as the prognostic variable or the transporting
velocity, and 3) transforming the form, typically by scaling
it. Filtering enables different types of terms to be separated and
handled by different schemes, as required by Algorithms \ref{alg:SIQN} and \ref{alg:splitphysics}. Detecting particular components
of the weak form enables their values to be updated, for example, the
transport schemes in Algorithm \ref{alg:SIQN} need to have the
transporting velocity updated at each iteration. Transforming forms
enables timestepping schemes to transfer the explicit terms to the
right hand side of an equation by scaling them by a multiple of $-dt$
where $dt$ is the model timestep. All of these manipulations require
us to be able to loop over the forms in the residual, query the
associated labels and perform actions based on their presence or their
values. Next we will describe the objects required to define labels and
filters before illustrating this with some simple examples.

The finite element discretisation of each term in an equation is
represented in Firedrake by a UFL form object. FML defines a
\texttt{Label} object which is used to tag a UFL form to create a
\texttt{Term}, with the intention that the \texttt{Term} represents
one or more forms in the residual. The \texttt{Label} object consists
of a name and an optional value. Labels can be used to tag terms to
identify them as representing particular processes, such as transport,
or they can identify particular components of the form. For example, the value of the \texttt{subject}
label is the \texttt{Function} in the \texttt{ufl.Form}
that represents the prognostic variable, or `subject' of the
equation, to be solved for. These operations are composable, so more than one label can
be attached to a term. The end result of labelling the terms in the
residual is a \texttt{LabelledForm} object that holds a list of terms,
pairing individual \texttt{ufl.Form} objects with FML labels. Labelled
forms can be added together, multiplied by floating point numbers and
extended by adding additional terms. More complex manipulations are
enabled by the \texttt{label\_map} method that provides a way of
filtering the terms in the form based on their labels, applying
different mappings to terms that pass or fail the conditions defined
in the filter. A common operation might be to filter for terms that
have a certain label, such as \texttt{transport}, drop terms that do
not have this label and multiply those that do by the timestep. We
shall demonstrate this in the following example.

\subsubsection{Example 1: IMEX-Euler applied to the 1d advection diffusion equation}

We demonstrate the use of FML by applying a split timestepping scheme
to the 1D advection diffusion equation
\begin{equation}\label{eq:1dadvdiff}
  q_t + cq_x - \kappa q_{xx} = 0,
\end{equation}
where $q$ is a scalar field, $c$ is the constant transporting velocity
and $\kappa$ is the constant diffusivity. The time discretisation
applies an explicit timestep of forward Euler to the transport term
$cq_x$ and an implicit timestep of backward Euler to the diffusion
term $\kappa q_{xx}$:
\begin{align}
  q^* &= q^n - c \Delta t q_x^n, \label{eq:explicit-advdiff} \\
  q^{n+1} - \kappa \Delta t q_{xx}^{n+1} &= q^*, \label{eq:implicit-advdiff}
\end{align}
where $q^n$ is the value of $q$ at time step $n$, $q^*$ is an
intermediate value for $q$ and $\Delta t$ is the timestep.

We illustrate the use of FML for this example in three steps: 1)
define and label the weak form of the residual, 2) build solvers for
the explicit and implicit systems, 3) run the timestepping loop. Since
the intention of this example is to demonstrate FML rather than the
compatible finite element discretisation, we use the space of piecewise
linear continuous Lagrange polynomials, denoted by $P_1$, for the
finite element discretisation as this greatly simplifies the
forms. The finite element discretisation of equation~(\ref{eq:1dadvdiff}) is then
\begin{equation}\label{eq:weak-advdiff}
  \int_\Omega vq_t \diff x + c \int_\Omega v q_x \diff x + \kappa \int_\Omega v_x q_x \diff x = 0, \quad \forall v \in P_1,
\end{equation}
where $v\in P_1$ is the test function and the diffusion term has been
integrated by parts.

Listing \ref{lst:define-residual} presents the code to define this
spatial discretisation of the residual with the $P_1$ function space
denoted by \texttt{V}. Here we omit the code that sets up the mesh,
function spaces and parameters $c$ and $\kappa$, although this can be
found in \citet{thomas_bendall_2026_24173215}. The first two lines define the
test function \texttt{v} and function \texttt{q}, corresponding to $v$
and $q$ in equation~(\ref{eq:weak-advdiff}). The next four lines define
the residual, i.e. the left side of equation~(\ref{eq:weak-advdiff}),
using the FML labels \texttt{time\_derivative}, \texttt{transport} and
\texttt{diffusion} to tag the corresponding terms. The final line
labels the function \texttt{q} as the subject of this residual,
i.e. the prognostic field. Note that the labels used for tagging just
take the UFL form as their argument and their value is simply the
default, \texttt{True}, in contrast to the \texttt{subject} label that
has the value \texttt{q}.

\begin{minipage}{0.94\linewidth}
\vspace{0.5em} 
\begin{lstlisting}[ language={[firedrake]python},
  style=framed, frame=single, firstline=24, lastline=31, firstnumber=24, caption={This
  code defines and labels the residual corresponding to
  equation~(\ref{eq:weak-advdiff}).  The full script for this example
  is available in \citet{thomas_bendall_2026_24173215}.}
]
q = Function(V)
v = TestFunction(V)
residual = (
    time_derivative(q * v * dx) +
    transport(c * q.dx(0) * v * dx) +
    diffusion(kappa * q.dx(0) * v.dx(0) * dx)
    )
residual = subject(residual, q)    # label q as the subject
\end{lstlisting}
\label{lst:define-residual}
\end{minipage}

To define the solvers for the explicit and implicit timesteps, equation~(\ref{eq:explicit-advdiff}) and equation~(\ref{eq:implicit-advdiff}),
we need to create UFL forms given by the following equations
\begin{align}
  \label{eq:explicit-advdiff_weak}
  \int_\Omega vq^* \diff x &= \int_\Omega vq^n \diff x - \Delta t c \int_\Omega v q_x^n \diff x, &\quad \forall v \in P_1, \\
  \label{eq:implicit-advdiff_weak}
  \int_\Omega vq^{n+1} \diff x  + \Delta t \kappa \int_\Omega v_x q_x^{n+1} \diff x &= \int_\Omega vq^* \diff x, &\quad \forall v \in P_1.
\end{align}
To do this we use the \texttt{label\_map} method to loop over the
terms in the residual making use of the predefined mappings
\texttt{drop} and \texttt{keep} to either remove or retain the term
based on the value of the filter. This code is given in Listings
\ref{lst:define-explicit} and \ref{lst:define-implicit}. In Listing
\ref{lst:define-explicit} we create the left and right sides of equation~(\ref{eq:explicit-advdiff_weak}). The left and right sides both contain
the time derivative, or mass, term so first we extract this (lines 2-4) by using the filter to check if each term has the
label \texttt{time\_derivative}, keeping the one that does and
dropping all others. For the left side of the
\texttt{LinearVariationalProblem}, we replace the subject of the mass term with a trial
function using the predefined \texttt{replace\_subject} function, a
lightweight wrapper around the UFL replace method which extracts the
subject of the term and replaces it with the new expression, in this
case simply a trial function \texttt{q\_trial} (lines 6-7). The right side of the
explicit solver consists of two forms, one coming from the discretisation of
the time derivative (the mass term) and the other from the transport term. Both are
evaluated using the values of $q$ at time step $n$. We first
extract the transport term and multiply it by $-\Delta t$ (lines 9-11), where the minus sign is required to move the term from
the left side to the right side. We then add this to the mass term we extracted earlier and replace the
subject with the function \texttt{qn} that holds the values of $q$ at
time step $n$. We can then extract the UFL forms from
\texttt{lhs\_explicit} and \texttt{rhs\_explicit} to create a
\texttt{LinearVariationalProblem} and corresponding
\texttt{LinearVariationalSolver} that, when solved, will put the
result into the function \texttt{qstar}.

\begin{minipage}{0.94\linewidth}
\vspace{0.5em}
\begin{lstlisting}[
   language={[firedrake]python},
   style=framed, frame=single,
   firstline=39, lastline=59, firstnumber=39,
   caption={This code defines the solver for the explicit transport
  step: \texttt{qstar = qn - dt*transport(qn)} i.e. equation~(\ref{eq:explicit-advdiff_weak}).  \texttt{qn} and \texttt{qstar}
  are Firedrake \texttt{Function}s with \texttt{qn} holding the values
  of $q$ at time step $n$ and \texttt{qstar} holding the values of $q$
  after the explicit solve.}
]
# extract the time derivative, or mass, term from the residual
mass = residual.label_map(lambda t: t.has_label(time_derivative),
                          map_if_true=keep,
                          map_if_false=drop)
# the lhs of the explicit problem is the mass term with unknown q_trial
lhs_explicit = mass.label_map(all_terms,
                              map_if_true=replace_subject(q_trial))
# extract the transport term from the residual and multiply by -dt
transport_rhs = residual.label_map(lambda t: t.has_label(transport),
                                   map_if_true=lambda t: -dt*t,
                                   map_if_false=drop)
# the rhs of the explicit problem is mass + transport_rhs evaluated at qn
rhs_explicit = mass + transport_rhs
rhs_explicit = rhs_explicit.label_map(all_terms,
                                      map_if_true=replace_subject(qn))
# setup the problem and solver, storing the result in qstar
prob_explicit = LinearVariationalProblem(lhs_explicit.form,
                                         rhs_explicit.form,
                                         qstar)
solver_explicit = LinearVariationalSolver(prob_explicit)
\end{lstlisting}
\label{lst:define-explicit}
\end{minipage}

The implicit solver is defined in a similar way (see listing
\ref{lst:define-implicit}). To create the left side of
equation~(\ref{eq:implicit-advdiff_weak}) we begin with the original
residual, extract the diffusion term and multiply it by $\Delta t$
(lines 2-4). We then add this to the mass term defined in listing
\ref{lst:define-explicit} and replace the subject with a trial
function (lines 6-8). The right side only contains the form coming
from the discretisation of the time derivative, i.e. the mass term
defined in listing \ref{lst:define-explicit}, evaluated using the
result of the explicit solve, \texttt{qstar} (lines 10-11). We now
create the \texttt{LinearVariationalProblem} and corresponding
\texttt{LinearVariationalSolver} that, when solved, will this time put
the result into the function \texttt{qnp1} (lines 13-16).

\begin{minipage}{0.94\linewidth}
\vspace{0.5em}
\begin{lstlisting}[
   language={[firedrake]python},
   style=framed, frame=single,
   firstline=60, lastline=75, firstnumber=60,
   caption={This code defines the solver for the implicit diffusion
  step: \texttt{qnp1 + dt*diffusion(qnp1) = qstar} i.e. equation~(\ref{eq:implicit-advdiff_weak}). \texttt{qstar} and \texttt{qnp1}
  are Firedrake \texttt{Function}s with \texttt{qstar} holding the
  values of $q$ after the explicit solve and \texttt{qnp1} holding the
  values of $q$ at after the implicit solve.}
]
# extract the diffusion term from the residual and multiply by dt
diffusion_lhs = residual.label_map(lambda t: t.has_label(diffusion),
                                   map_if_true=lambda t: dt*t,
                                   map_if_false=drop)
# the lhs of the implicit problem is mass + diffusion_lhs with unknown q_trial
lhs_implicit = mass + diffusion_lhs
lhs_implicit = lhs_implicit.label_map(all_terms,
                                      map_if_true=replace_subject(q_trial))
# the rhs of the implicit problem is the mass term evaluated at qstar
rhs_implicit = mass.label_map(all_terms,
                              map_if_true=replace_subject(qstar))
# setup the problem and solver, storing the result in qnp1
prob_implicit = LinearVariationalProblem(lhs_implicit.form,
                                         rhs_implicit.form,
                                         qnp1)
solver_implicit = LinearVariationalSolver(prob_implicit)
\end{lstlisting}
\label{lst:define-implicit}
\end{minipage}

The final part of the code is the time loop which calls the solve
method of first the explicit solver and then the implicit solver, see
Listing \ref{lst:timeloop}. In practice, the code in listings
\ref{lst:define-explicit}-\ref{lst:define-implicit} would be in a
function that could be applied to any residual labelled with the
appropriate labels.

\begin{minipage}{0.94\linewidth}
\vspace{0.5em}
\begin{lstlisting}[
   language={[firedrake]python},
   style=framed, frame=single,
   firstline=80, lastline=84, firstnumber=80,
   caption={This code is the timeloop implementing the scheme defined in equations~(\ref{eq:explicit-advdiff_weak})-(\ref{eq:implicit-advdiff_weak}) using
   the solvers defined in listings \ref{lst:define-explicit}-\ref{lst:define-implicit}.}
]
while t < tmax:
    t += dt
    solver_explicit.solve()
    solver_implicit.solve()
    qn.assign(qnp1)
\end{lstlisting}
\label{lst:timeloop}
\end{minipage}

This simple example demonstrates the power of manipulating the terms
in the equation based on the presence or absence of certain
labels. The full strength of FML becomes even more apparent when we
exploit information stored in the label value. We use this in two
important situations in Gusto. The first situation is when we need to
evaluate source terms that come from physics schemes. The second
situation is when we require a linearisation of the residual, for
example when applying a Quasi-Newton method to solve the full
nonlinear equations or when testing a new algorithm on the linear
equations. We illustrate both situations in the following two
examples.

\subsubsection{Example 2: Evaluation of complex source terms} \label{sec:physics_demo}

This example demonstrates the use of FML for encapsulating
instructions for the evaluation of complex source terms in the value
of the label. This is typically how source terms from physics
parametrisation schemes are implemented in Gusto. Each scheme is
implemented in its own Python class which, at instantiation, adds the
required terms to the equation residual form. The class contains an
\texttt{evaluate} method that is then used by the time discretisation to
update the physics source term during the timestepping loop. This means that the same physics scheme
can be used with different equation sets (where appropriate) and time discretisations. The
simplest example of this functionality is the addition of a source
term to the transport equation for specific humidity of water vapour $q_v$

\begin{equation}
  \frac{\partial q_v}{\partial t} + \boldsymbol{u}\cdot\nabla q_v - \mathcal{S}(q_v, \boldsymbol{u}, q_{sat}) = 0,
\end{equation}
where the source term $\mathcal{S}$ is a function of the water vapour
$q_v$, the wind velocity $\boldsymbol{u}$ and the saturation specific humidity
$q_{sat}$, which can itself also be a function of the atmospheric
conditions. $\mathcal{S}$ models the source of water vapour due to
evaporation and is given by

\begin{equation}
  \mathcal{S} =
  \begin{cases}
    C(q_{sat} - q_v), & \textrm{if } \, q_{sat} > q_v, \\
    0, & \textrm{otherwise},
  \end{cases}
\end{equation}
where $C$ is scaling factor that may be dependent on atmospheric
properties such as the surface wind velocity and the bulk transfer
coefficient. This gives a source of water vapour if the current water
vapour is less than the saturation specific humidity. For simplicity,
in the following we shall assume that both $q_{sat}$ and $C$ are
constant values specified when the \texttt{Evaporation} class is
instantiated. Note that in this case the source term does not contain
any derivatives and therefore its weak form can be defined by
multiplying by a test function from the appropriate space and
integrating over the domain.

In Listing \ref{lst:physics_demo} we outline how the source term
$\mathcal{S}$ is defined, added to the residual and then
evaluated. The source term is defined in the \texttt{\_\_init\_\_}
method of the class. This definition is purely symbolic: the \texttt{Function}s
\texttt{self.source} and \texttt{self.q} have not yet been allocated
any values. The \texttt{evaluate} method assigns these values when it
is called at the appropriate point of the timestepping scheme. We see
this in Listing \ref{lst:physics_demo} where the \texttt{ForwardEuler}
timestepping class extracts the physics schemes (in this case only the
evaporation scheme) by checking for the presence of a \texttt{physics}
label and, if found, adding its value, i.e. its \texttt{evaluate}
method, to a list. When the forward Euler \texttt{apply} method is
called, the first action is to update the source terms by calling the
stored evaluation methods, passing in the current model state
\texttt{x\_in} and the timestep (not used in this example).

\begin{minipage}{0.94\linewidth}
\vspace{0.5em}
\begin{lstlisting}[style=framed, frame=single, language={[firedrake]python},
  caption={Sketch of the use of FML for evaluating source terms. The
  \texttt{Evaporation} physics class defines the components required
  for the weak form of the source term and its \texttt{evaluate}
  method implements the instructions for updating the source term
  given the model fields \texttt{x} and the timestep \texttt{dt}. The
  update happens at the appropriate point in the timestep, as
  illustrated below by the \texttt{ForwardEuler} class.}
  ]
  class Evaporation(PhysicsParametrisation}:

      def __init__(self, equation, qsat, C):
          # identify the component of the mixed function that holds the water vapour field
          self.q_index = equation.field_names.index("water_vapour")
          # get the function space for the water vapour from the mixed function space
          Vq = equation.function_space.sub(self.q_index)

          v = TestFunction(Vq)
          self.source = Function(Vq)  # this is a placeholder for the source term S
          self.q = Function(Vq)       # this is a placeholder for the water vapour field

          # the symbolic expression to evaluate the source term
          self.source_expr = C * conditional(qsat > self.q, (qsat - self.q), 0)

          # add the source term to the equation residual, labelling it as a physics term with an
          # associated evaluate method
          equation.residual -= physics(v * self.source * dx, self.evaluate)
    
      def evaluate(self, x, dt):
          # update the values in self.q to be those from the coming state x
          self.q.assign(x.subfunctions[self.q_index])
          # interpolate the source expression to the source function use in the term
          self.source.interpolate(self.source_expr)


  class ForwardEuler(TimeDiscretisation):

      def __init__(self, equation...):
          ...
          # this list will store the methods required to evaluate the source terms
          self.evaluate_source = []
      
          # check for any terms in the residual that have the physics label
          for term in equation.residual:
              if term.has_label(physics):
                  # get the value of the physics label, i.e. the evaluate method, and
                  # append it to the list of sources to evaluate 
                  self.evaluate_source.append(term.labels[physics])

      def apply(x_out, x_in):
          # this method applies one timestep of the forward Euler method to x_in, storing the
          # result as x_out.
          for evaluate in self.evaluate_source:
              # evaluate the source term, given the current model state and timestep
              evaluate(x_in, dt)
          # continue the timestep with the updated values of the source terms
          ...
\end{lstlisting}
\label{lst:physics_demo}
\end{minipage}

This example can be easily extended to evaluate more complex source
terms. Including the wind dependence in $C$ requires the velocity
field to be extracted from the equation and added to the symbolic expression on line 14 and the values
of the velocity to be updated in the \texttt{evaluate} method, as is
done for the water vapour. To allow the saturation specific humidity
to depend on the atmospheric conditions, the \texttt{qsat} argument
could become a Python function that returns $q_{sat}$ and is called
from the \texttt{evaluate} method, passing the model state \texttt{x}.

\subsubsection{Example 3: Generating linearised forms} \label{sec:linearisation demo}

Linearised forms are an essential component of many solver
algorithms. For example, the Semi-Implicit Quasi-Newton timestepping
scheme described in Section \ref{sec:design} requires an approximate
Jacobian generated from a linearisation of the residual (the $S$ in
line 8 of Algorithm \ref{alg:SIQN}, and equation~(\ref{eqn:linear_system})). Linearisation is also used for generating
linear equation sets for testing algorithms.

The example in Listing \ref{lst:linearisation_demo} shows how we
generate the linearisation of each term in the residual, store it as
the value of the \texttt{linearisation} label for that term and later
retrieve it using the \texttt{label\_map} method. The linearisation
can be generated using UFL by differentiating the form corresponding
to the term with respect to the prognostic variable \texttt{X} (the
`subject'). The use of FML enables some additional functionality
beyond using this default linearisation. Firstly, we can specify the
state that we wish to linearise about and update this at a specified
frequency, rather than always using the current model state. This is
more efficient because the matrix factorisations for the solver do not
have to be updated at every timestep. We do this by replacing any
instances of \texttt{X} remaining after the linearisation by
\texttt{X\_ref}, the reference profile that we wish to linearise
about, as shown in Listing \ref{lst:linearisation_demo}. Secondly, we
can specify which forms to include in the linearisation based on their
labels, making it straightforward to completely drop transport terms,
for example. Thirdly, FML provides a mechanism for entirely overriding the UFL-generated linearisation by using a different linear form as the value of the
\texttt{linearisation} label, rather than just the same form linearised around a different state as shown above. This is
useful when the form is complex and retaining the entire linearisation
is inefficient, but the term describes a physical process that is important to capture in the linear solve. For example, including just the vertical component of
buoyancy terms in the linear solve of the SIQN algorithm produces a
solver that has the same convergence properties as that with the full
linearisation but is more efficient \citep{melvin2024mixed}. This
combination of functionality enabled by FML makes Gusto a flexible and
extensible tool for investigating linear solvers.

\begin{minipage}{0.94\linewidth}
\vspace{0.5em}
\begin{lstlisting}[style=framed, frame=single, language={[firedrake]python}, language=python, caption={Sketch of the use of FML for generating, storing and
   retrieving the linearisation of the residual.}]
  X = Function(V)               # where V is the function space for this problem
  X_ref = Function(V)           # to store the profiles we wish to linearise about
  residual = subject(..., X)    # definition of the residual as a labelled form with subject X
  for term in residual:
      # generate the linearisation of the form and retain the original labels
      linear_term = Term(ufl.derivative(term.form, X), term.labels)
      # replace X with X_ref
      linear_term_with_ref = replace_subject(X_ref)(linear_term)
      # store as the value of the linearisation label for the original term
      term = linearisation(term, linear_term_with_ref)
      ...
# extract the linearisation of the residual
  linearised_residual = residual.label_map(
      lambda t: t.has_label(linearisation)
      map_if_true=lambda t: Term(t.get(linearisation).form, t.labels),
      map_if_false=drop)
\end{lstlisting}
\label{lst:linearisation_demo}
\end{minipage}

\section{Results}
\label{sec:results}
The aim of this section is to present a selection of results that
showcase the tools available in Gusto, such as the choice of grids,
$H(\textrm{div})$ element families, order of finite element spaces, different
spatial discretisation options and timestepping schemes.  As these
options can be combined in so many ways, it is only possible to
present results from a tiny fraction of the configurations.  We have
chosen to use a selection from the hierarchy of standard test cases
that span the stages of development of dynamical cores.  Although we
do include some results here that demonstrate the numerical properties
of certain configurations, the purpose of this section is not to
provide a rigorous verification of any particular choice of
discretisation. Such results are demonstrated elsewhere; for instance
\citet{cotter2023compatible2} describes the spatial discretisation of
the compressible Euler equations in a vertical slice and
\citet{bendall2020compatible} describes the moist discretisation, both
used in Section \ref{sec:compressible Euler results}.  The settings
used for each result from this section are detailed in Appendix
\ref{sec:test_case_settings}, but initial conditions are generally not
described where these followed a standard test case from the literature.

\subsection{Terminator Toy chemistry}
\begin{figure}[htp!]
\centering
\includegraphics[width=0.9\textwidth]{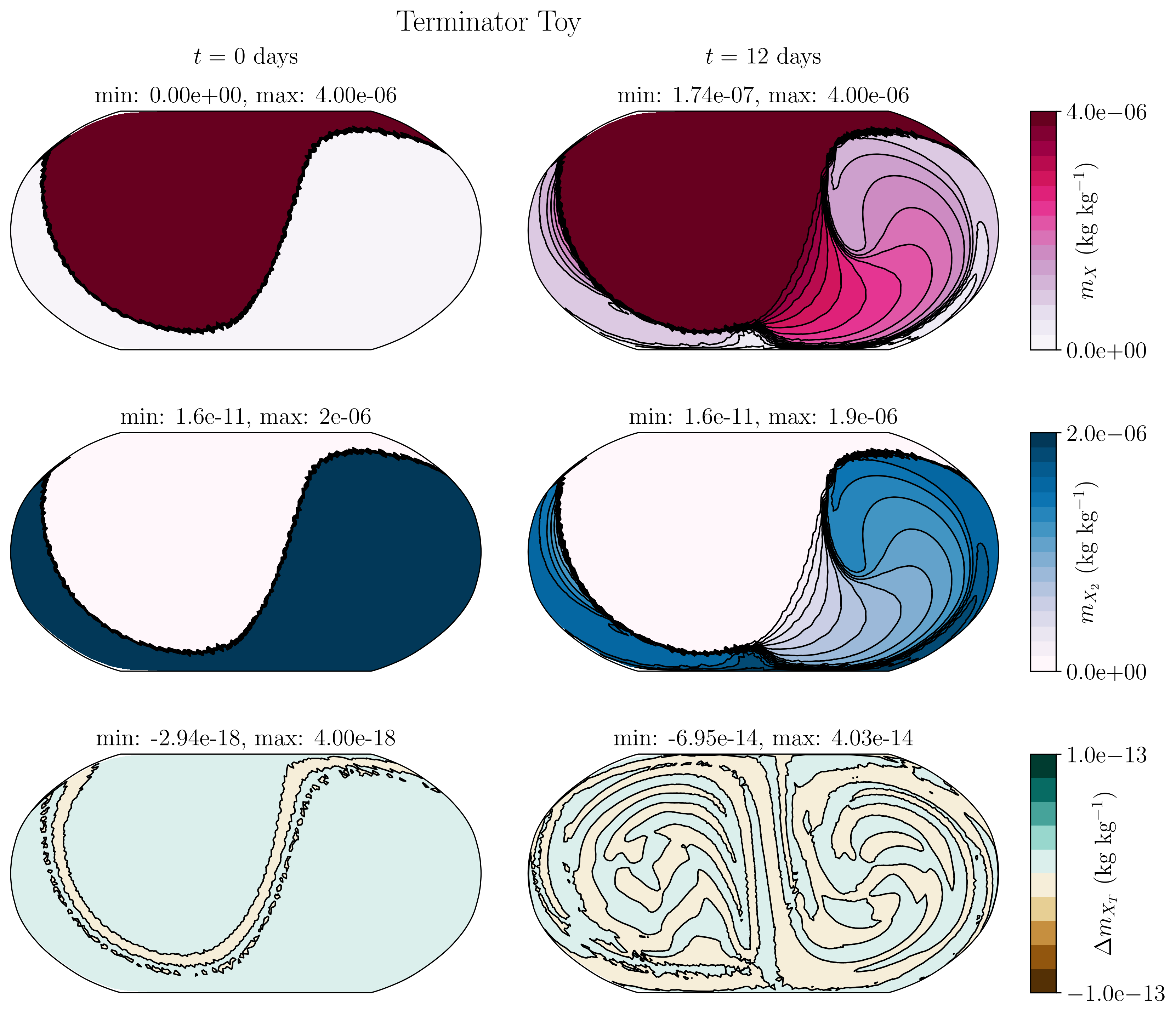}
\caption{The evolution of the chemical species from the ``Terminator Toy'' test of \citet{lauritzen2015terminator}.
The fields at the initial time are shown on the left, and at the final time ($t=12$ days) which follows a single revolution around the sphere by a deformational flow.
The bottom row shows the change in the total $m_T=m_X+2m_{X_2}$ from the initial time, which remains very close to zero.
The mixing ratio $m_X$ is contoured every $2.5\times 10^{-7}$ kg kg$^{-1}$, $m_{X_2}$ is contoured every $1.25\times 10^{-7}$ kg kg$^{-1}$ while
$m_T$ is contoured every $2\times 10^{-14}$.
}
\label{fig:terminator}
\end{figure}

As described in Section \ref{sec:transport}, transport schemes are an essential component in the solution of the equations given in Section \ref{sec:equations} and it is useful to be able to test this part of the software independently from the other terms in the equations. Here we demonstrate a transport scheme in Gusto using
 the ``Terminator Toy'' test of \citet{lauritzen2015terminator}, which involves the transport of two coupled tracers $X$ and $X_2$ on the surface of a sphere.
The chemicals interact such that the sum $m_T=m_X+2m_{X_2}$ should be preserved, where $m_X$ and $m_{X_2}$ are the mixing ratios of the individual chemicals.
When positivity or monotonicity limiters are applied to the numerical transport scheme, this may break the correlation between the two species which will be highlighted by this test.
This test is also a demonstration of simple ``physics-dynamics coupling'', as the conversion term between the species can be handled separately from the transport terms. \\
\\
For this test, we use the DG upwind scheme described in Section \ref{subsec:DGtransport} to transport the species on an icosahedral sphere mesh, and apply the vertex-based limiter of \citet{kuzmin2010vertex}.
The term describing the chemical reaction is separated from the transport term using a simple sequential splitting, as described in Section \ref{sec:physics_demo}. \\
\\
Figure \ref{fig:terminator} shows the mixing ratios of the two individual species at the initial time and following a single rotation around the sphere by a deformational flow.
The change in the sum of the species is also displayed, with differences from zero remaining close to machine precision, which demonstrates that this
transport scheme does preserve the correlation between the tracers.

\subsection{Shallow water}
\begin{figure}[htp!]
\centering
\includegraphics[width=0.95\textwidth]{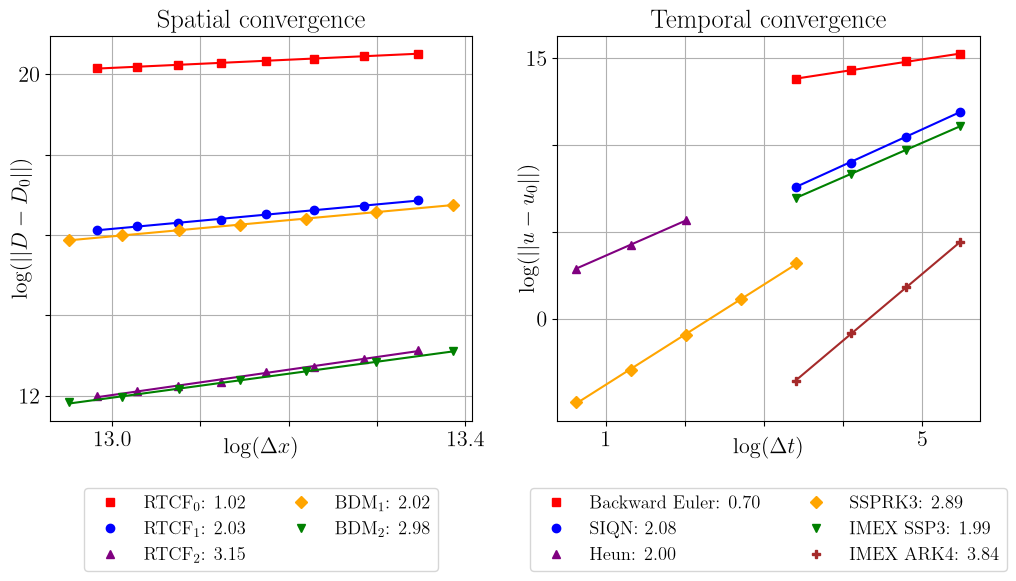}
\caption{Convergence results using the shallow water equations.
(Left) the spatial convergence of a number of orders and families of finite element, all using the same timestepper.
The test case used is a variant of the steady geostrophic flow of \citet{williamson1992standard},
with a perturbation added to initiate gravity waves.
(Right) the temporal convergence of different timesteppers, all using the same spatial discretisation,
measured with the Rossby-Haurwitz wave test case of \citet{williamson1992standard}.
Explicit timestepping schemes use a smaller timestep to ensure that they remain stable.
For both sets of convergence results, the legend includes the gradient of the line of best fit, which
indicates the order of accuracy. Measured convergence rates are close to or match the theoretical orders.}
\label{fig:sw_convergence}
\end{figure}

\begin{figure}[htp!]
\centering
\includegraphics[width=0.95\textwidth]{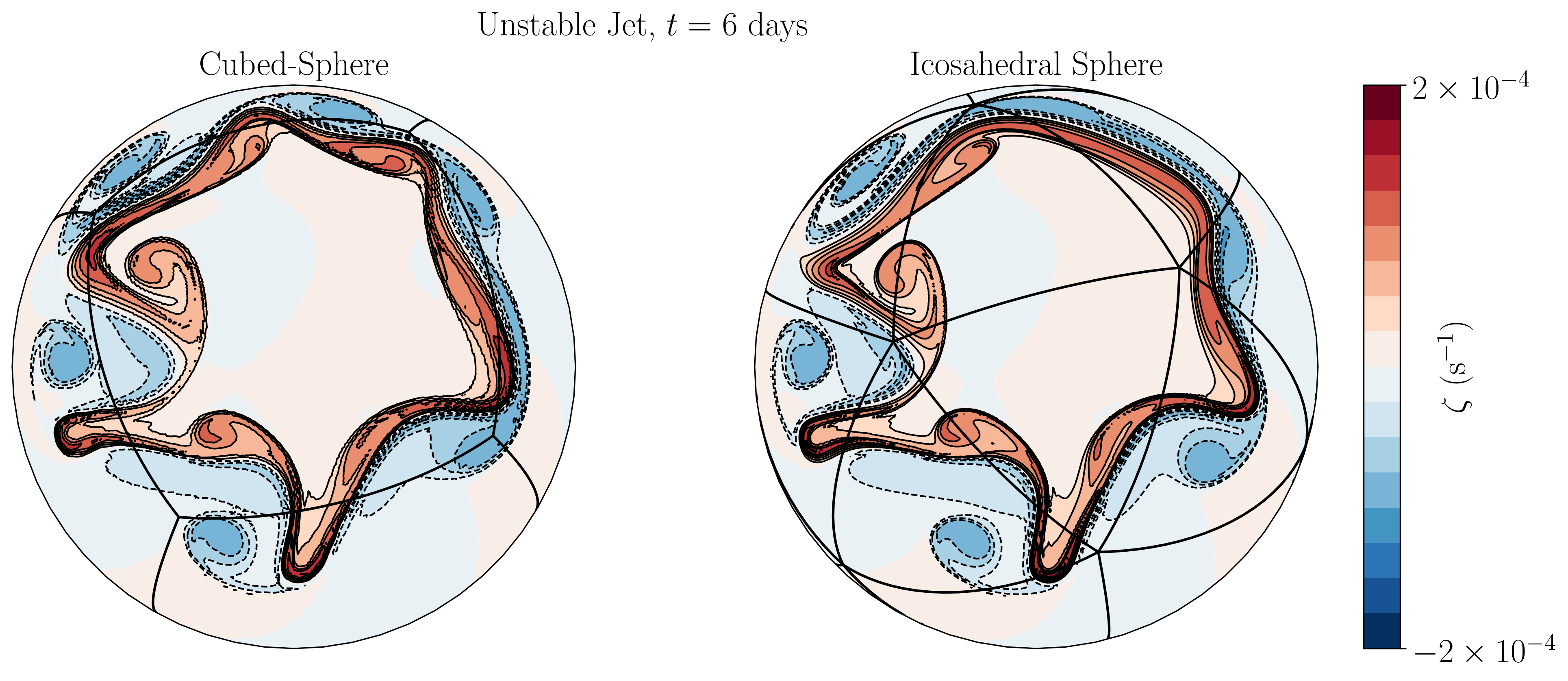}
\caption{A view of the relative vorticity around the North pole at day 6 of a simulation of the Galewsky jet test case.
(Left) the RT1 elements on a cubed-sphere mesh. (Right) the BDM1 elements on an icosahedral-sphere mesh.
Both cases use the Semi-Implicit Quasi-Newton timestepper. The edges of the panels of the cube and icosahedron are highlighted,
and the relative vorticity is contoured every $1.25\times 10^{-5}$ s$^{-1}$. As is common for this test case, the zero contour line
is omitted to avoid emphasising small variations around the reference state.}
\label{fig:galewsky}
\end{figure}

The standard suite of shallow water test cases originally proposed in
\citet{williamson1992standard} provides a common starting point for
testing new discretisations. This test suite is often augmented by
the unstable barotropic jet described in
\citet{galewsky2004initial}. \\
\\
Here we present two sets of convergence results (spatial and temporal)
for a variety of configurations measured using shallow water equations.
Both spatial and temporal convergence tests involve computing a high resolution reference solution (in space or time)
for the depth or velocity field, against which errors are measured using the $L^2$-norm. \\
\\
Figure \ref{fig:sw_convergence} shows the spatial convergence of the depth field
for a number of choices of finite element for the prognostic variables,
using a Trapezium Rule discretisation (also known as Crank-Nicolson).
The spatial convergence was measured using the linear shallow-water equations, with a variation of the
geostrophically-balanced flow test of \citet{williamson1992standard} in which a perturbation is added
to trigger gravity waves. This test is described in Appendix \ref{sec:test_case_settings} as the ``Gravity Wave'' test.
The details of the initial conditions are given in Appendix \ref{sec:spatial_convergence_ics}. \\
\\
Temporal convergence for a range of timestepping schemes are also shown in Figure \ref{fig:sw_convergence},
using the first order RT elements on quadrilateral cells and measuring the error in the velocity field.
Note that the Semi-Implicit Quasi-Newton scheme is iterated to convergence, by using 8 iterations for its outer loop
and a single iteration for each inner loop.
The temporal convergence test used the Rossby-Haurwitz wave test case of \citet{williamson1992standard}, with the nonlinear
shallow water equations.
The errors converge at or close to the theoretical rates for different finite element families and timestepping schemes.
\\
 \\
We show results from the unstable jet of \citet{galewsky2004initial} in Figure \ref{fig:galewsky}, using the nonlinear shallow water equations.
These results use the Semi-Implicit Quasi-Newton timestepper described in Section \ref{sec:siqn},
and compare the first order RT elements on a cubed-sphere mesh with the first-order BDM elements on an icosahedral-sphere. \\
Although much of the jet structure is very similar between the two meshes, there are some visible
differences that can be associated with the mesh structure.

\subsection{Compressible Euler equations} \label{sec:compressible Euler results}
\begin{figure}[htp!]
\centering
\includegraphics[width=0.6\textwidth]{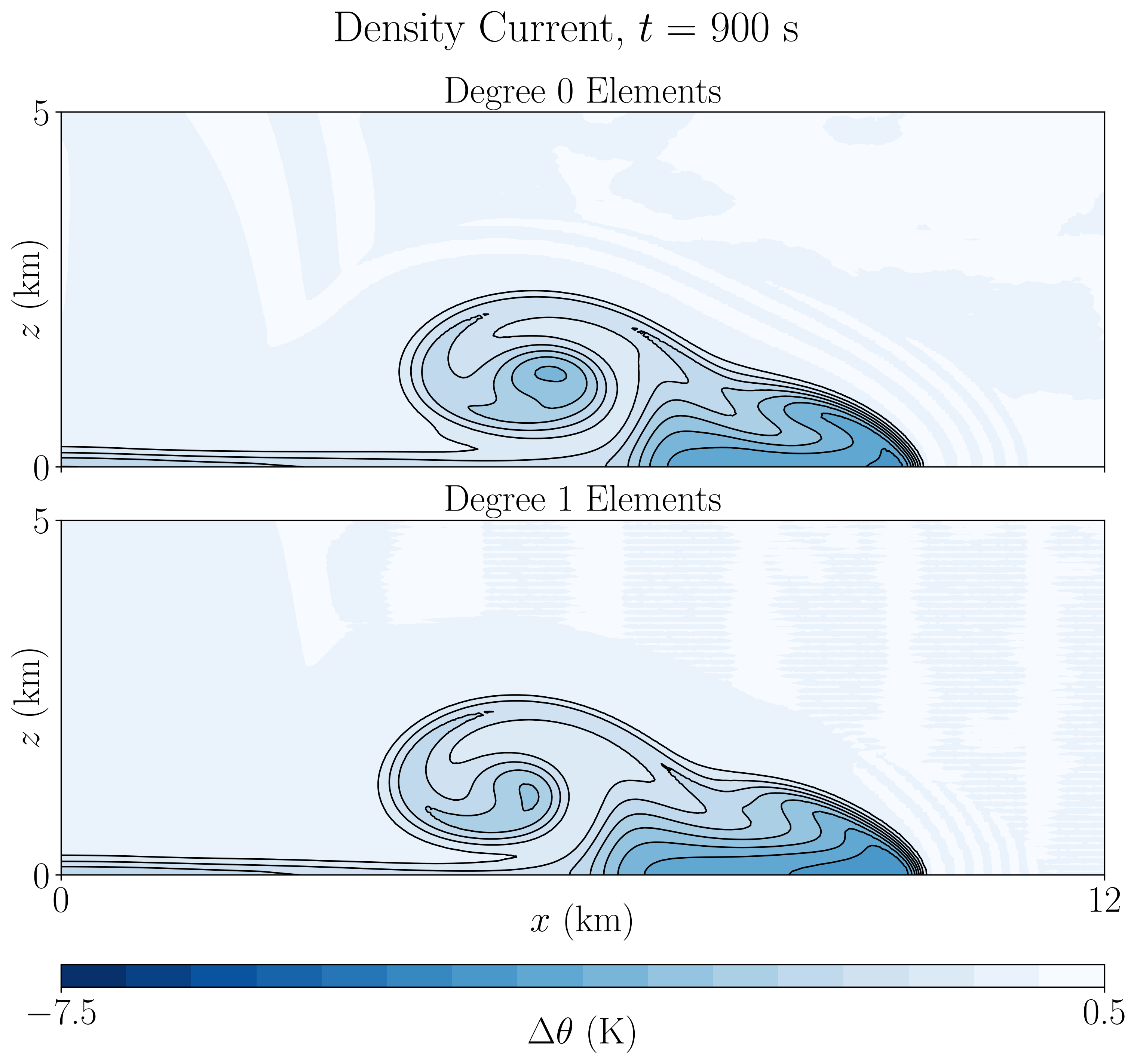}
\caption{The perturbation in potential temperature field after 900 s from the falling density current test of \citet{straka1993numerical}.
(Top) with the lowest-order finite elements, and (bottom) with the next-to-lowest-order finite elements.
The potential temperature is contoured every $0.5$ K, and as is common for this test case,
the zero contour line omitted to avoid emphasising small variations around the reference state.
}
\label{fig:straka}
\end{figure}

\begin{figure}[htp!]
\centering
\includegraphics[width=0.95\textwidth]{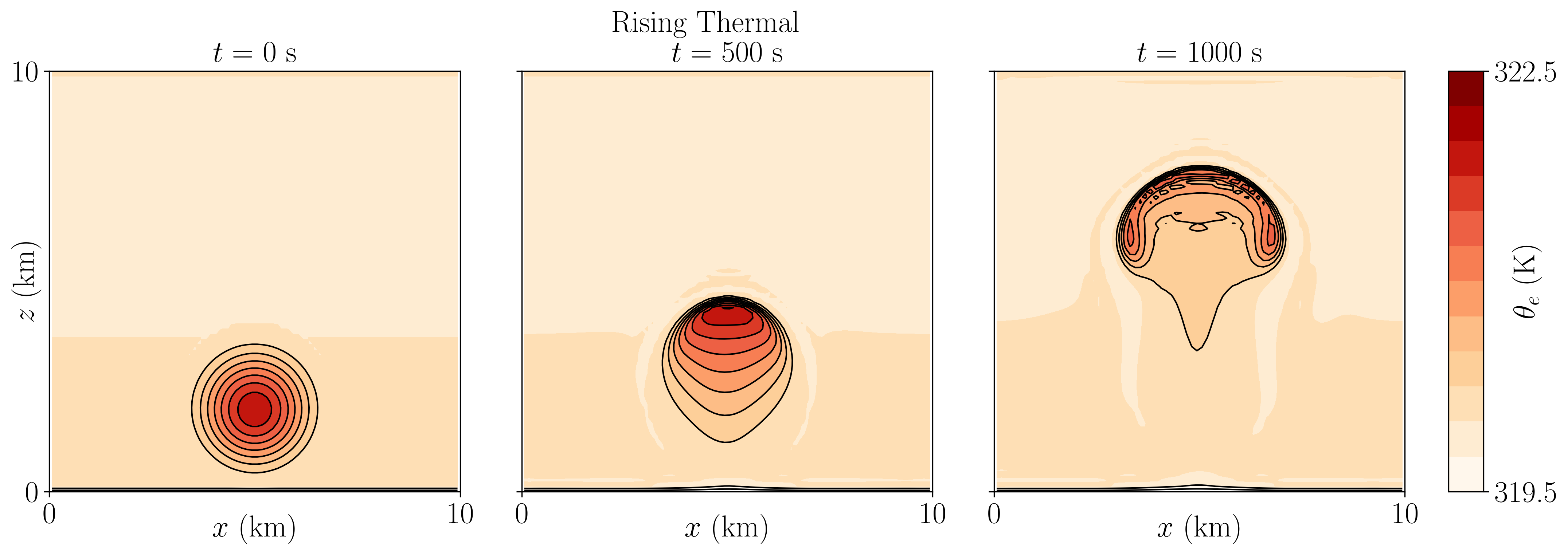}
\caption{Vertical cross-section of the equivalent potential temperature field from the three-dimensional moist rising thermal test,
at a series of times. The equivalent potential temperature is contoured every $0.5$ K,
with the $320$ K contour line omitted to avoid emphasising small variations around the reference state.}
\label{fig:3d_bubble}
\end{figure}

\begin{figure}[htp!]
\centering
\includegraphics[width=0.95\textwidth]{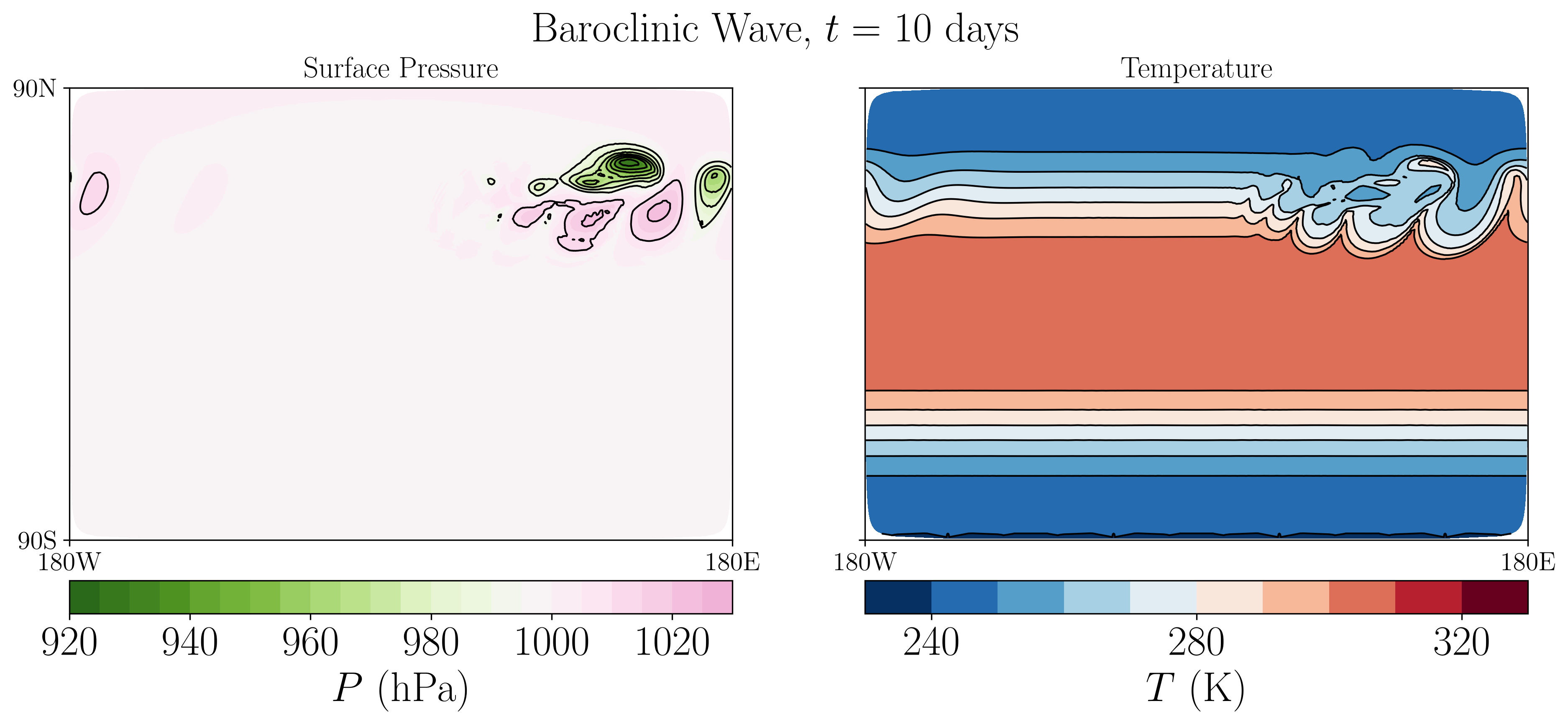}
\caption{The surface pressure and temperature after 10 days of a deep baroclinic wave simulation on the cubed-sphere.
The surface pressure is contoured with colours every 5 hPa and lines every 10 hPa,
and the 1000 hPa contour line is omitted to avoid emphasising small variations around the reference state. The temperature
is contoured every 10 K.}
\label{fig:baroclinic}
\end{figure}

We have also implemented many of the standard test cases for the non-hydrostatic compressible Euler equations that are used in modern weather and climate models.
In this section we show results from a small selection of these, which demonstrate different domains in which these methods work.
Throughout this section we use the Semi-Implicit Quasi-Newton timestepper, with the settings
described in Appendix \ref{sec:siqn_details}.
 \\
\\
The first test is the falling density current of \citet{straka1993numerical} in a vertical slice. This test involves the addition of diffusion to both the potential temperature and velocity equations.
Figure \ref{fig:straka} shows the final potential temperature field for simulations using the lowest-order and next-to-lowest order finite elements on the vertical slice.
The lowest-order methods require the transport and diffusion terms to use the ``recovered'' scheme of \citet{bendall2019recovered} to attain second-order accuracy.
\\
\\
To demonstrate the use of moist physics parametrisations, we use a moist three-dimensional rising thermal simulation in an extruded doubly-periodic plane.
The initial conditions are those of \citet{bendall2020compatible}, which is an extension of the dry case of \citet{kelly2012continuous} and \citet{melvin2019mixed}.
A cross-section of the final equivalent potential temperature field is shown in Figure \ref{fig:3d_bubble}.
As in the two-dimensional test of \citet{bryan2002benchmark}, the atmosphere is saturated and cloudy everywhere, so that the moist processes are active in the whole domain.
A simple saturation adjustment parametrisation is used, with the details described in \citet{bendall2020compatible}.
\\
\\
The next test case is the deep baroclinic wave of \citet{ullrich2014proposed} on the sphere.
An initial state in thermal wind balance is perturbed to produce a mid-latitude baroclinic instability, providing a flow that is typical of those found in numerical weather prediction.
The surface pressure and temperature fields after 10 days are shown in Figure \ref{fig:baroclinic}, which compare well with those in \citet{melvin2024mixed} and demonstrate that our methods can replicate results from the literature.

\begin{figure}[htp!]
\centering
\includegraphics[width=0.95\textwidth]{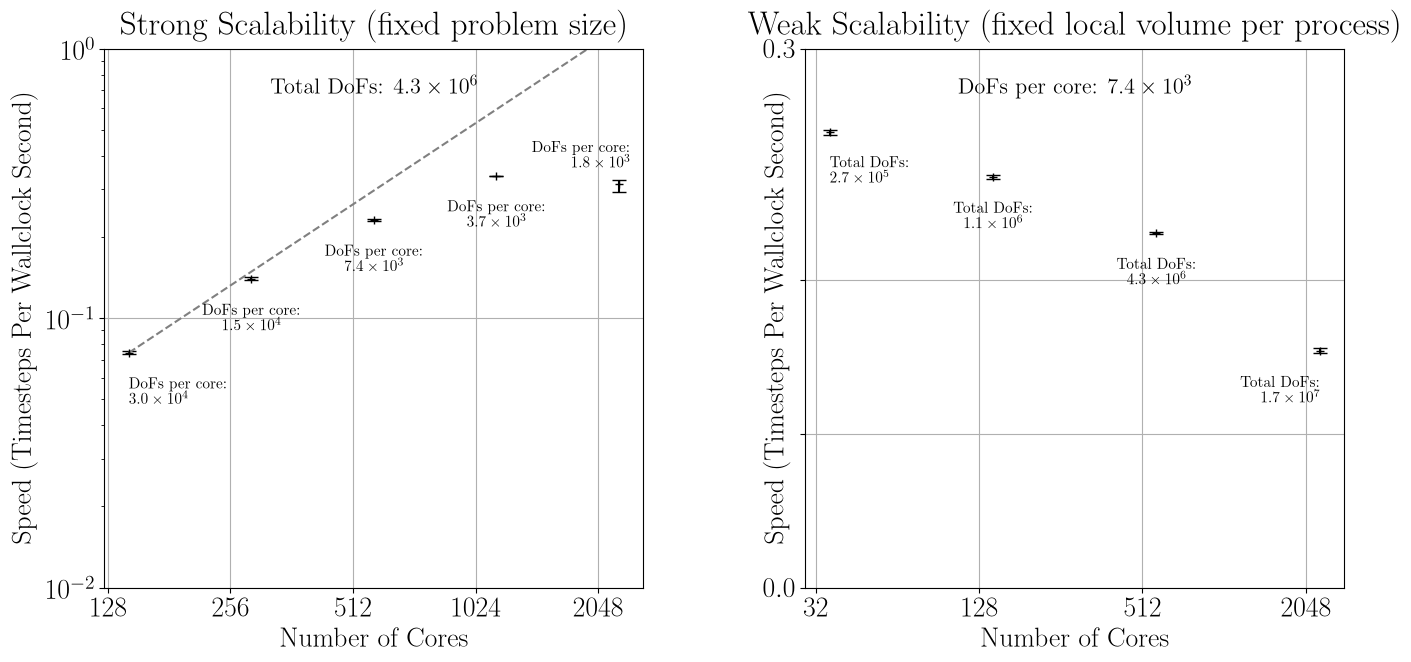}
\caption{Computational performance of Gusto's 3D compressible Euler model, demonstrating
scalability via computational speed (the inverse of wallclock time per timestep), measured
using the simulation of a fixed number of timesteps of the baroclinic wave configuration.
Full details of the configurations used for testing scalability are given in
Appendix \ref{sec:scalability_details}.
At each processor count, three measurements were performed and the range in values are shown
through the error bars in the figure.
Each point is labelled with the number of prognostic degrees of freedom (the sum of
the degrees of freedom for the $\boldsymbol{u}$, $\rho$ and $\theta$ fields), either
globally or per core, as appropriate.
(Left) strong scalability, as problem size is fixed at the C96 resolution and
computational resource is varied. The dashed line shows ideal scaling from the
first data point.
(Right) weak scalability, as the global problem size is varied with the computational
resource, to keep the local problem size constant. The number of layers was kept constant,
while the time step was varied to keep the horizontal Courant number constant
(the number of timesteps was fixed). The number of cores used per node was also constant.
Unlike the strong scaling plot, the $y$-axis does not use a logarithmic scale,
to aid interpreting drop-off in performance as the total number of cores increases.
}
\label{fig:scalability}
\end{figure}

To demonstrate the computational performance of Gusto, we show some strong and weak scalability results in Figure \ref{fig:scalability}.
These results also use the deep baroclinic wave test of \citet{ullrich2014proposed}, with a range of resolutions and different
computational resources.
Full details of the testing can be found in Appendix \ref{sec:scalability_details}.
This shows that Gusto can strong scale to thousands of cores,
inheriting the parallel scalability properties of Firedrake,
as demonstrated for example by \citet{betteridge2021code}.
We reach the limit of strong scalability at local DoF counts of around 2000,
where communication costs dominate.

\begin{figure}[htp!]
\centering
\includegraphics[width=0.95\textwidth]{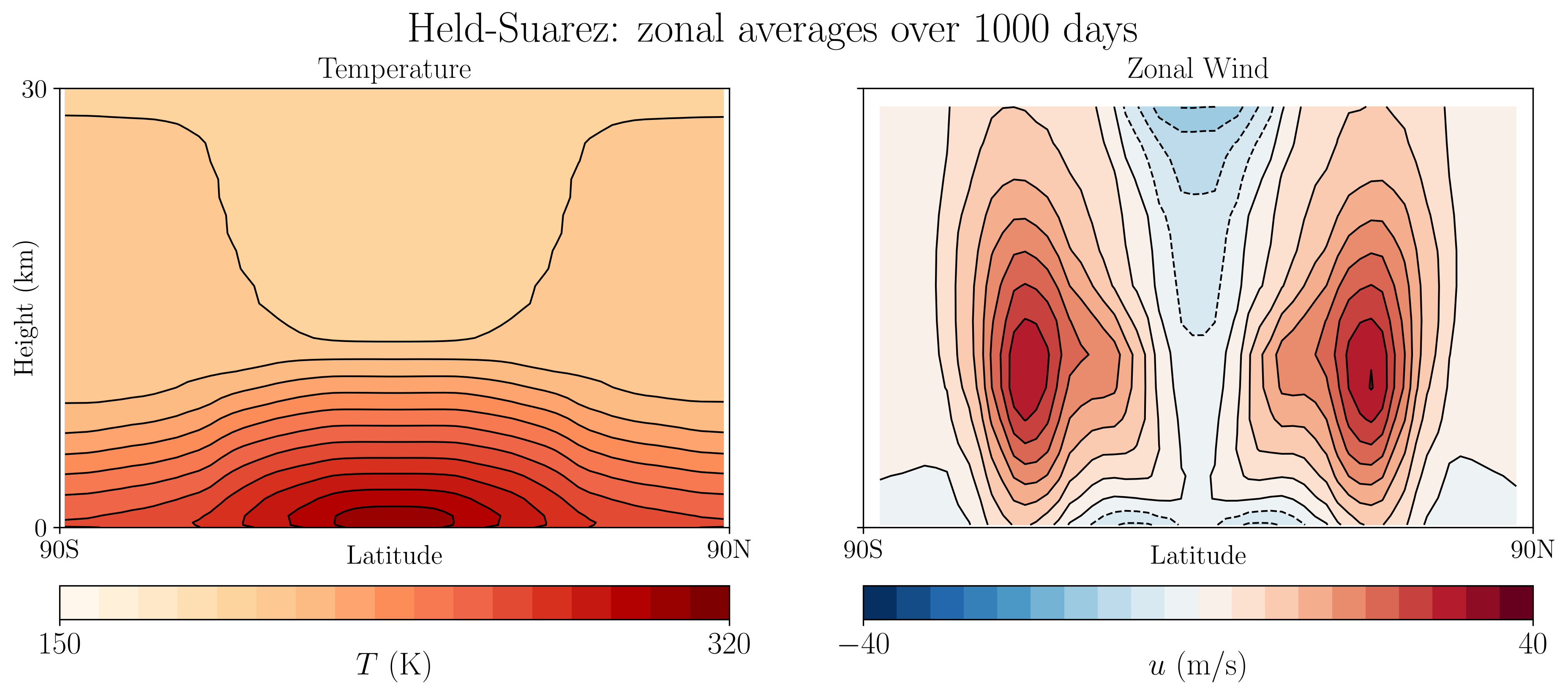}
\caption{Cross-sections of the temperature and zonal wind fields from the Held-
Suarez test case, averaged zonally and over 1000 days of simulation.
The temperature is contoured every 10 K and the wind is contoured every 4 m s$^{-1}$.}
\label{fig:held_suarez}
\end{figure}

Figure \ref{fig:held_suarez} shows results from the Held-Suarez test case.
This is a climate-like simulation with long time integration, running for 1200 days,
and including simple wind drag and temperature relaxation physics parametrisations.
The zonally-averaged temperature and zonal wind fields are shown in Figure \ref{fig:held_suarez},
with time-averaging applied over final 1000 days to exclude the spin-up period.
The Held-Suarez test uses the same initial conditions as the baroclinic wave.
The results in Figure \ref{fig:held_suarez} are similar to those from the literature
such as \citet{brown2024physics}, and demonstrate that Gusto can be run stably
for long time integrations.

\section{Summary}
\label{sec:summary}
This paper presents Gusto, a toolkit for developing and testing numerical
methods for use in weather and climate models.
Gusto uses compatible finite element methods, which are being used in the Met Office's next dynamical core due to their desirable numerical properties, and ability to facilitate
non-orthogonal meshes.
The purpose of Gusto is to provide a flexible framework in which new numerical
methods can be easily prototyped and explored across a range of equations and test problems.
In this paper, we have described how Gusto's modular code design facilitates this flexibility.
The Form Manipulation Language allows terms in equations to be labelled according to
how the terms will be discretised numerically.
A range of numerical results have been included to showcase Gusto's capabilities. \\
\\
Gusto is underpinned by the Firedrake software, which provides not only a wide range of finite elements but access to the PETSc solver library. Firedrake provides Gusto with a user-friendly interface to highly efficient and scalable implementation of finite element methods.
Gusto is also able to exploit other features of Firedrake, such as auto-differentiability and adjoint generation, which can be used to answer important contemporary research questions relating to data assimilation and
the fusion of numerical methods with data driven approaches. \\
\\



\subsection*{Code and data availability}{The Gusto website \url{https://www.firedrakeproject.org/gusto/} describes the project and includes documentation and examples. The specific version of Gusto used in the paper is available here: \url{https://doi.org/10.5281/zenodo.23419270}  \citep{thomas_bendall_2026_23419271}. The main Gusto code repository is \url{https://github.com/firedrakeproject/gusto}. Additional case studies can be found in a separate repository \url{https://github.com/firedrakeproject/gusto_case_studies}. For the specific components of the Firedrake project used in this paper, see \url{https://zenodo.org/doi/10.5281/zenodo.23957393}. For the input files of all examples and benchmarks presented, see \url{https://doi.org/10.5281/zenodo.24173214} \citep{thomas_bendall_2026_24173215} or the GitHub repository \url{https://github.com/jshipton/gusto_gmd_paper_results}. All plots are produced using the tomplot library, available at \url{https://github.com/tommbendall/tomplot}.} 





\appendix
\section{Test Case Settings}  \label{sec:test_case_settings}  
Here we summarise the configuration details for the test cases used to produce the results in Section \ref{sec:results}.
These details can be found in Table \ref{tab:test_case_settings}.

\begin{sidewaystable*}[tbp!]
\caption{Details of the configurations used for the test cases in Section \ref{sec:results}.}
\label{tab:test_case_settings}
\footnotesize
\begin{tabular}{l|ccccccc}
\tophline
Test Case & Equation Set & Timestepper & Finite Elements & \begin{tabular}{@{}c@{}}Horizontal \\ Mesh\end{tabular} & \begin{tabular}{@{}c@{}}Number \\ of  \ Layers\end{tabular} & Time step & Final Time \\
\middlehline
Terminator Toy & \begin{tabular}{@{}c@{}}Tracer transport \\ with chemistry\end{tabular} & Sequential Splitting & BDM1 & i16 & None & 450 s & 12 days \\
Gravity Wave & Linear Shallow Water & Trapezium Rule & Multiple & Multiple & None & 60 s & 1 hour \\
Gravity Wave (Reference) & Linear Shallow Water & Trapezium Rule & RT2 & C36 & None & 15 s & 1 hour \\
Rossby-Haurwitz Wave & Shallow Water & Multiple & RT1 & C128 & None & Multiple & 1 day \\
Rossby-Haurwitz Wave (Reference) & Shallow Water & RK4 & RT1 & C128 & None & 1 s & 1 day \\
Unstable Jet & Shallow Water & SIQN & BDM1 and RT1 & i64 and C100 & None & 150 s & 6 days \\
Density Current & Compressible Euler & SIQN & RT0 and RT1 & 100 m and 200 m & 64 and 32 & 1 s & 900 s \\
Rising Thermal & Moist Compressible Euler & SIQN & RT1 & $\Delta x = \Delta y =$ 0.125 km & 80 & 2 s & 1000 s \\
Baroclinic Wave & Compressible Euler & SIQN & RT1 & C48 & 15 & 600 s & 10 days \\
Held-Suarez & Compressible Euler & SIQN & RT1 & C12 & 15 & 1200 s & 1200 days \\
\bottomhline
\end{tabular}
\caption{Notation: RT$x$ refers to the de Rham complex that uses the Raviart-Thomas elements on quadrilateral cells for the $H(\textrm{div})$ function space, and results in an $x$-th order polynomials for the $L_2$ space. Similarly, BDM$x$ refers to the Brezzi-Douglas-Marini de Rham complex on triangular cells.
SIQN denotes the Semi-Implicit Quasi-Newton timestepper described in Section \ref{sec:siqn}. When describing horizontal meshes of the sphere, i$y$ and C$y$ correspond to icosahedral-sphere and cubed-sphere meshes with $y$ cells per edge of the icosahedron or cube panel.}
\end{sidewaystable*}

\section{Details of Scalability Tests}  \label{sec:scalability_details}  
The scalability tests of Section \ref{sec:compressible Euler results} used the deep
baroclinic wave test of \citet{ullrich2014proposed}, using the SIQN timestepper and RT1
elements with 15 vertical layers.
All tests were run on the Isambard 3 supercomputer using compute nodes equipped with two NVIDIA Grace CPUs (144 cores per node) connected by HPE Slingshot-11 interconnect. The number of cores per node was kept constant at 36 for all tests, which
ensures the highest resolution configurations fit within memory.
Speeds were calculated as the inverse of the time taken to complete one timestep,
with the first timestep excluded to avoid initialisation costs.
Measurements were repeated three times for each configuration.
\\
\\
Strong scaling tests used a C96 cubed-sphere mesh, with a timestep of 225 s and
a final time of 6 hours (96 timesteps).
Results were recorded with 144, 288, 576, 1152 and 2304 cores. \\
\\
Weak scaling tests used a range of horizontal resolutions, with the number of
vertical layers kept constant at 15, and the timestep varied to keep the
horizontal Courant number constant.
The local volume was fixed at 96 columns per core, and the number of timesteps
was kept constant at 96.
The timesteps, resolutions and core counts are recorded in Table \ref{tab:weak_scalability}. \\
\\
\begin{table*}[tbp!]
\caption{Details of the configurations used for the weak scalability tests.}
\label{tab:weak_scalability}
\footnotesize
\begin{tabular}{l|ccccccc}
\tophline
Horizontal Mesh & Time step & Final Time & Number of Cores \\
\middlehline
C24 & 900 s & 24 hours & 36 \\
C48 & 450 s & 12 hours & 144 \\
C96 & 225 s & 6 hours & 576 \\
C192 & 112.5 s & 3 hours & 2304 \\
\bottomhline
\end{tabular}
\end{table*}

\section{Initial Conditions for Shallow Water Gravity Wave} \label{sec:spatial_convergence_ics}
The spatial convergence test included in Section \ref{sec:results} uses a variation
of the geostrophically-balanced flow test of \citet{williamson1992standard}, in which
a perturbation is added to trigger gravity waves.
We use the linear shallow water equations for this test.
While the mean state is set to $\overline{\boldsymbol{u}}=\boldsymbol{0}$ and $\overline{D}=2000$ m,
the initial perturbation fields (to be solved for) are set to:
\begin{equation}
u'=u_0\cos\phi, \qquad v'=0, \qquad
D' = -\frac{r\Omega u_0}{g}\cos^2\phi + D_0e ^{-(\gamma/R)^2},
\end{equation}
where $u'$ and $v'$ are the zonal and meridional components of the velocity,
$\phi$ is the latitude and $\gamma$ is the great arc distance in radians from
longitude $\lambda=-\pi/2$ rad and latitude $\phi=\pi/3$ rad.
The perturbation has magnitude $D_0=100$ m and width $R=\pi/9$ rad.
The remaining parameters take the values of the geostrophically-balanced test
of \citet{williamson1992standard}: the planet's radius is $r=6.37122\times10^6$ m,
the planet's rotation rate is given by $\Omega=7.292\times10^{-5}$ s$^{-1}$,
and acceleration due to gravity is set to $g=9.80616$ m s$^{-2}$.

\section{Details of Semi-Implicit Quasi-Newton Timestepper} \label{sec:siqn_details}
This Appendix describes some important details of the SIQN
timestepper, which follow the formulation of the Gungho dynamical
core.  These options are all described in detail elsewhere, but are
summarised here to aid reproducibility.  Each feature is implemented
as an option within Gusto, and all those described here have been used
in the numerical tests in Section \ref{sec:results}.  The
Quasi-Newton structure uses 2 outer and 2 inner loops by default, with the
off-centering parameter set to $\alpha=1/2$.  The inner loop involves
the solution of a linear problem, equation~(\ref{eqn:linear_system}), of finite
element discretisations of the form:
\begin{subequations}
\begin{align}
\boldsymbol{u}'+\alpha\Delta t \mathcal{F}(\chi')=\boldsymbol{u}_r, \\
q'+\tau\Delta \mathcal{T}_{\boldsymbol{u}'}(q')=q_r, \label{eqn:linear_solver_scalar}
\end{align}
\end{subequations}
where equation~(\ref{eqn:linear_solver_scalar}) is the generic form for a scalar prognostic $q$,
such as depth, density or potential temperature.
Equation~(\ref{eqn:linear_solver_scalar}) includes the relaxation parameter $\tau$ rather than the off-centering parameter $\alpha$,
and as in \citet{melvin2024mixed}, the relaxation parameter $\tau$ is set to 1.
\\
\\
For all-but-the-first inner iteration (in practice, the second inner iteration of each outer iteration),
any non-velocity residuals (e.g. for the depth, density or potential temperature fields) are set to zero.
The residuals for these fields would otherwise not be updated within the inner loop, since
these equations do not include forcing terms, and setting these residuals to zero improves
the convergence of the inner loop.
This approach is described in \citet{melvin2024mixed}.
Related to this is the choice to set the relaxation parameter $\tau=1$, which follows
\citet{wood2014inherently} and \citet{melvin2024mixed}. \\
\\
The linear problem used in the inner loop requires the definition of reference profiles.
For the numerical results shown in this paper, the reference profiles are set to the values
of the field from the previous timestep, as in \citet{melvin2024mixed}. \\
\\
A key conclusion of \citet{wood2026replicating} is that it is possible to replicate
the desirable dispersion properties of a Semi-Implicit Semi-Lagrangian (SISL) scheme
whilst using Eulerian transport schemes via the use of a ``density predictor''.
This involves replacing the transported depth or density field with
\begin{equation}
\rho^n \to \rho^n - \frac{1}{2}\Delta t \rho^n \nabla\cdot\boldsymbol{u}^n.
\end{equation} \\
\\
The transport of all variables is timestepped using the Strong Stability Preserving
3-Stage Runge-Kutta (SSPRK3) scheme \citep{shu1988efficient}, and any fields in
conservative form (such as depth or density) use the ``advective-then-flux''
formulation described in \citet{bendall2023solution} and \citet{melvin2024mixed},
in which only the final stage of the SSPRK3 scheme uses the conservative form of
the transport equation.








\subsection*{Author Contribution}Conceptualisation: JS, TMB, CJC, DAH. Funding acquisition: JS, CJC. Investigation: JS, TMB, CJC. Methodology: JS, TMB, CJC. Project administration: JS, TMB, CJC. Software: JS, TMB, DMA, TCA, AB, JB, CJC, THG, DAH, JH-C, NH, JAP, CJW, DW. Supervision: JS, TMB, CJC. Validation: JS, TMB, TCA, AB, NH, DW. Visualisation: TMB, AB. Writing - original draft: JS, TMB. Writing - review and editing: JS, TMB, CJC, DMA, DAH, JH-C.



\subsection*{Acknowledgements}
This research and software development has been supported by UK
Research and Innovation (UKRI) through grants from the Engineering and
Physical Sciences Research Council (EP/R029628/1, EP/W026066/1, EP/W029731/1), the Natural
Environment Research Council (NE/I02013X/1, NE/K006789/1 and
NE/R008795/1), the Science and Technology Facilities Council (UKRI/ST/B000495/1) and the UK Exascale computing programme ExCALIBUR
EX20-8.

The authors acknowledge the use of resources provided by the Isambard
3 Tier-2 HPC Facility. Isambard 3 is hosted by the University of
Bristol and operated by the GW4 Alliance (https://gw4.ac.uk) and is
funded by UK Research and Innovation; and the Engineering and Physical
Sciences Research Council (EP/X039137/1). The authors are grateful to
the entire Firedrake development team for their support and advice
throughout the development of Gusto.







\bibliography{references.bib}

\end{document}